\input harvmac 
\input epsf.tex

\overfullrule=0mm
\def\file#1{#1}
\def\figbox#1#2{\epsfxsize=#1\vcenter{
\epsfbox{\file{#2}}}} 
\newcount\figno
\figno=0
\def\fig#1#2#3{
\par\begingroup\parindent=0pt\leftskip=1cm\rightskip=1cm\parindent=0pt
\baselineskip=11pt
\global\advance\figno by 1
\midinsert
\epsfxsize=#3
\centerline{\epsfbox{#2}}
\vskip 12pt
{\bf Fig.\the\figno:} #1\par
\endinsert\endgroup\par
}
\def\figlabel#1{\xdef#1{\the\figno}}
\def\encadremath#1{\vbox{\hrule\hbox{\vrule\kern8pt\vbox{\kern8pt
\hbox{$\displaystyle #1$}\kern8pt}
\kern8pt\vrule}\hrule}}


\def\IR{\relax{\rm I\kern-.18em R}}
\font\cmss=cmss10 \font\cmsss=cmss10 at 7pt

\font\cmss=cmss10 \font\cmsss=cmss10 at 7pt
\def\IZ{\relax\ifmmode\mathchoice
{\hbox{\cmss Z\kern-.4em Z}}{\hbox{\cmss Z\kern-.4em Z}}
{\lower.9pt\hbox{\cmsss Z\kern-.4em Z}}
{\lower1.2pt\hbox{\cmsss Z\kern-.4em Z}}\else{\cmss Z\kern-.4em Z}\fi}
\def\IN{\relax{\rm I\kern-.18em N}}
\def\b{\circ}
\def\n{\bullet}

\def\Rww{R^{\b \b}}

\def\Rwb{R^{\b \n}}

\def\gbbbb{\Gamma_4^{\hbox{$\scriptstyle \b \b$}\kern -8.2pt
\raise -4pt \hbox{$\scriptstyle \b \b$}}}
\def\gnnnn{\Gamma_4^{\hbox{$\scriptstyle \n \n$}\kern -8.2pt  
\raise -4pt \hbox{$\scriptstyle \n \n$}}}
\def\gnnnnnn{\Gamma_6^{\hbox{$\scriptstyle \n \n \n$}\kern -12.3pt
\raise -4pt \hbox{$\scriptstyle \n \n \n$}}}
\def\gbbbbbb{\Gamma_6^{\hbox{$\scriptstyle \b \b \b$}\kern -12.3pt
\raise -4pt \hbox{$\scriptstyle \b \b \b$}}}
\def\gbbbbc{\Gamma_{4\, c}^{\hbox{$\scriptstyle \b \b$}\kern -8.2pt
\raise -4pt \hbox{$\scriptstyle \b \b$}}}
\def\gnnnnc{\Gamma_{4\, c}^{\hbox{$\scriptstyle \n \n$}\kern -8.2pt
\raise -4pt \hbox{$\scriptstyle \n \n$}}}
\def\Rbud#1{{\cal R}_{#1}^{-\kern-1.5pt\blacktriangleright}}
\def\Rleaf#1{{\cal R}_{#1}^{-\kern-1.5pt\vartriangleright}}
\def\Rbudb#1{{\cal R}_{#1}^{\circ\kern-1.5pt-\kern-1.5pt\blacktriangleright}}
\def\Rleafb#1{{\cal R}_{#1}^{\circ\kern-1.5pt-\kern-1.5pt\vartriangleright}}
\def\Rbudn#1{{\cal R}_{#1}^{\bullet\kern-1.5pt-\kern-1.5pt\blacktriangleright}}
\def\Rleafn#1{{\cal R}_{#1}^{\bullet\kern-1.5pt-\kern-1.5pt\vartriangleright}}
\def\Wleaf#1{{\cal W}_{#1}^{-\kern-1.5pt\vartriangleright}}
\def\Cleaf{{\cal C}^{-\kern-1.5pt\vartriangleright}}
\def\Cbud{{\cal C}^{-\kern-1.5pt\blacktriangleright}}
\def\Crleaf{{\cal C}^{-\kern-1.5pt\circledR}}


\Title{\vbox{\hsize=3.truecm \hbox{SPhT/03-148}}}
{{\vbox {
\bigskip
\centerline{Geodesic Distance in Planar Graphs:}
\medskip
\centerline{An Integrable Approach}
}}}
\bigskip
\centerline{ 
P. Di Francesco\foot{philippe@spht.saclay.cea.fr}}
\medskip
\centerline{ \it Service de Physique Th\'eorique, CEA/DSM/SPhT}
\centerline{ \it Unit\'e de recherche associ\'ee au CNRS}
\centerline{ \it CEA/Saclay}
\centerline{ \it 91191 Gif sur Yvette Cedex, France}
\bigskip
\vskip .5in
\noindent 
We discuss the enumeration of planar graphs using bijections with suitably decorated trees,
which allow for keeping track of the geodesic distances between faces of the graph. The
corresponding generating functions obey non-linear recursion relations on the geodesic
distance. These are solved by use of stationary multi-soliton tau-functions of suitable reductions
of the KP hierarchy. We obtain a unified formulation of the (multi-) critical continuum limit
describing large graphs with marked points at large geodesic distances,
and obtain integrable differential equations for the corresponding scaling functions.
This provides a continuum formulation of two-dimensional quantum gravity, in terms 
of the geodesic distance.
\Date{10/03}

\nref\BIPZ{E. Br\'ezin, C. Itzykson, G. Parisi and J.-B. Zuber, {\it Planar
Diagrams}, Comm. Math. Phys. {\bf 59} (1978) 35-51.}
\nref\DGZ{P. Di Francesco, P. Ginsparg
and J. Zinn--Justin, {\it 2D Gravity and Random Matrices},
Physics Reports {\bf 254} (1995) 1-131.}
\nref\EY{B. Eynard, {\it Random Matrices}, Saclay Lecture Notes (2000),
available at {\sl http://www-spht.cea.fr/lectures\_notes.shtml} }
\nref\TUT{W. Tutte, 
{\it A Census of planar triangulations}
Canad. Jour. of Math. {\bf 14} (1962) 21-38;
{\it A Census of Hamiltonian polygons}
Canad. Jour. of Math. {\bf 14} (1962) 402-417;
{\it A Census of slicings}
Canad. Jour. of Math. {\bf 14} (1962) 708-722;
{\it A Census of Planar Maps}, Canad. Jour. of Math. 
{\bf 15} (1963) 249-271.}
\nref\KPZ{V.G. Knizhnik, A.M. Polyakov and A.B. Zamolodchikov, {\it Fractal Structure of
2D Quantum Gravity}, Mod. Phys. Lett.
{\bf A3} (1988) 819-826; F. David, {\it Conformal Field Theories Coupled to 2D Gravity in the
Conformal Gauge}, Mod. Phys. Lett. {\bf A3} (1988) 1651-1656; J.
Distler and H. Kawai, {\it Conformal Field Theory and 2D Quantum Gravity},
Nucl. Phys. {\bf B321} (1989) 509-527.}
\nref\SCH{G. Schaeffer, {\it Bijective census and random 
generation of Eulerian planar maps}, Electronic
Journal of Combinatorics, vol. {\bf 4} (1997) R20; see also
G. Schaeffer, {\it Conjugaison d'arbres
et cartes combinatoires al\'eatoires} PhD Thesis, Universit\'e 
Bordeaux I (1998).}
\nref\CENSUS{J. Bouttier, P. Di Francesco and E. Guitter, {\it Census of planar
maps: from the one-matrix model solution to a combinatorial proof},
Nucl. Phys. {\bf B645}[PM] (2002) 477-499.}
\nref\BMS{M. Bousquet-M\'elou and G. Schaeffer,
{\it Enumeration of planar constellations}, Adv. in Applied Math.,
{\bf 24} (2000) 337-368.}
\nref\BFG{J. Bouttier, P. Di Francesco and E. Guitter,
{\it Counting colored Random Triangulations}, 
Nucl.Phys. {\bf B641} (2002) 519-532.}
\nref\PS{D. Poulalhon and G. Schaeffer, 
{\it A note on bipartite Eulerian planar maps}, preprint (2002),
available at {\sl http://www.loria.fr/$\sim$schaeffe/}}
\nref\NOUSHARD{J. Bouttier, P. Di Francesco and E. Guitter, {\it
Combinatorics of hard particles on planar maps},
Nucl. Phys. {\bf B655} (2003) 313-341.}
\nref\CONC{M. Bousquet-M\'elou and G. Schaeffer, {\it The degree distribution
in bipartite planar maps: application to the Ising model}, preprint
math.CO/0211070.}
\nref\CS{P. Chassaing and G. Schaeffer, {\it Random Planar Lattices and 
Integrated SuperBrownian Excursion}, preprint (2002), to appear in 
Probability Theory and Related Fields, math.CO/0205226.}
\nref\TWOWALL{J. Bouttier, P. Di Francesco and E. Guitter, {\it  Random trees between two walls:
Exact partition
function}, J. Bouttier, P. Di Francesco and E. Guitter, Saclay preprint t03/086
and cond-mat/0306602 (2003), to appear in J. Phys. A: Math. Gen. (2003).}
\nref\GEOD{J. Bouttier, P. Di Francesco and E. Guitter, {\it Geodesic distance in planar graphs},
Nucl. Phys. {\bf B 663[FS]} (2003) 535-567.}
\nref\KKMW{H. Kawai, N. Kawamoto, T. Mogami and Y. Watabiki, {\it Transfer Matrix 
Formalism for Two-Dimensional Quantum Gravity and Fractal Structures of Space-time}, 
Phys. Lett. B {\bf 306} (1993) 19-26.}
\nref\AW{J. Ambj\o rn and Y. Watabiki, {\it Scaling in quantum gravity},
Nucl.Phys. {\bf B445} (1995) 129-144.}
\nref\AJW{J. Ambj\o rn, J. Jurkiewicz and Y. Watabiki, 
{\it On the fractal structure of two-dimensional quantum gravity},
Nucl.Phys. {\bf B454} (1995) 313-342.}
\nref\JM{M. Jimbo and T. Miwa, {\it Solitons and infinite dimensional Lie
algebras}, Publ. RIMS, Kyoto Univ. {\bf 19} No. 3 (1983) 943-1001, 
eq.(2.12).}
\nref\GD{I. Gelfand and L. Dikii, {\it Fractional powers of operators and
Hamiltonian systems}, Funct. Anal. Appl. {\bf 10:4} (1976) 13.}
\nref\EK{B. Eynard and C. Kristjansen, {\it Exact Solution of the O(n) 
Model on a Random Lattice}, Nucl.Phys. {\bf B455} (1995) 577-618, and {\it
More on the exact solution of the O(n) model on a random lattice 
and an investigation of the case $|n|>2$}, Nucl.Phys. {\bf B466} (1996) 463-487.}

\newsec{Introduction}

The work presented in this note was initially motivated by the need to better understand
on a combinatorial level the various results obtained via matrix models on the 
enumeration of graphs with fixed topology \BIPZ\ (see also \DGZ\ and \EY\
and references therein). The early work on this subject dates back to
a combinatorist, W. Tutte \TUT, who managed to enumerate many of the planar versions of these  
using recursion relations in the spirit of what we call today ``loop equations" for
matrix models. Such an approach, though combinatorial, failed to really explain
the simplicity of the algebraic equations determining the generating functions
for planar graphs.
Another motivation comes from the physics of two-dimensional quantum gravity. At the discrete
level, coupling matter to gravity simply amounts to define a statistical model (typically
with local Boltzmann weights) on a fluctuating base space,
in the form of random discretized surfaces. The continuum version of this involves
field-theoretical
descriptions of random surfaces with critical matter \KPZ, via the coupling of conformal
field theories (matter) to the Liouville field theory (metrics of the underlying space).

The interpretation of the planar graph results remained elusive
until the groundbreaking work of G. Schaeffer \SCH, who finally gave a beautifully
simple combinatorial explanation for these algebraic equations, in terms of decorated
trees. The idea was simply to establish bijections between classes of planar graphs
and suitably decorated trees, then easily enumerated via algebraic relations
obeyed by their generating functions. This technique proved quite general, 
and was extended to many classes of planar graphs, including graphs of
arbitrary valence \CENSUS, special classes of bipartite graphs called constellations
\BMS\ \BFG\ and other classes of bipartite graphs \PS, including the particular cases
of hard objects on planar graphs \NOUSHARD\ and of the Ising model on planar graphs \CONC. 

The great advantage of this bijective enumeration is that it allows for keeping
track of some details of the graphs in the language of trees. An important example of this
concerns the geodesic distance between say the vertices of random quadrangulations,
namely planar graphs with only tetravalent faces. In \CS, it was shown that the
geodesic distance of all vertices from an origin vertex of the graph may translate into
integer vertex labels in some corresponding trees, themselves realizing a discrete
version of the Brownian snake, further studied and extended
in \TWOWALL\ in the language of spatial branching processes. 
In \GEOD, it was shown that the generating functions
for planar graphs with two external legs obeyed non-linear recursion relations
on the maximal geodesic distance between the legs, and it is the 
purpose of this note to clarify and extend the results
of this paper.  

Note that no continuum (field-theoretical) treatment of two-dimensional
quantum gravity in terms of geodesic distance is available yet, only some partial
results were obtained using a transfer matrix formalism \KKMW\ \AW\ \AJW\ leading to some
conjectural continuum expression for a scaling function of the geodesic distance in the
case of pure gravity without matter. The present work provides an alternative
solution and gives access to a host of scaling functions for various (multi-) critical
models of matter systems coupled to two-dimensional quantum gravity.

The paper is organized as follows. In Sects. 2,3 and 4, we recall some known facts
on the bijective enumeration of planar graphs with respectively even valences,
arbitrary (even and odd) valences, and bicolored vertices. 
This relies on an iterative cutting procedure which, starting from a planar graph,
produces a decorated rooted tree. In turn the tree may be closed back in a unique way into
a planar graph, and we use this bijection to recover the algebraic relations
satisfied by the various counting functions involved. 
Sect. 5 is devoted to the introduction of the geodesic distance in the various
enumeration problems at hand: we present a unified picture involving formal
operators $Q$ generating the descendent trees around a vertex, and allowing
for keeping track of the geodesic distance from the root to the external face of the planar
graph to which the tree closes back.
The main results are recursion relations on the geodesic distance satisfied
by the generating functions for planar graphs with legs.  
In Sect. 6, we solve exactly a number of these recursion relations, by 
expanding the generating functions at large maximal geodesic distance $n$,
and resumming the resulting series. These display a remarkable ``integrable" structure
in that they generically involve tau-functions of the KP hierarchy. 
Sect. 7 is devoted to the critical continuum limit of the problem, in which
we consider large planar graphs and large geodesic distances. 
The generating functions calculated in Sect. 6 are shown to yield universal
scaling functions, characteristic of the various (multi-) critical points
of random surfaces with matter and with marked points at a fixed geodesic distance. 
We propose a generalization of our results based on differential equations obeyed by
the scaling functions, and illustrate it in the case
of the Ising model on random surfaces.
Finally, we gather a few concluding remarks in Sect. 8, where in particular
we discuss the integrable structure of our recursion relations, and of their
continuum counterparts.

\newsec{Planar graphs and trees I: the case of even valences}

In this section, we recall some results on the enumeration of planar graphs of even valence
with two extra ``legs". This is done via a bijection between planar graphs and
decorated trees, easily enumerated. 

\subsec{Tetravalent graphs}

\fig{The bijection between planar tetravalent graphs with two legs (a) 
and rooted ternary blossom-trees (d) is obtained via the following cutting procedure.
We first visit all edges of the graph in counterclockwise direction (b) and cut them iff
the resulting graph is still connected: we have indicated the succession of cut edges
by their number in the order of visit, from 1 to 6 here. The cut edges are then replaced by
pairs of white and black leaves (c). Moreover the incoming leg is replaced
by a white leaf and the outcoming one by a root, finally yielding a rooted ternary
blossom-tree (d).}{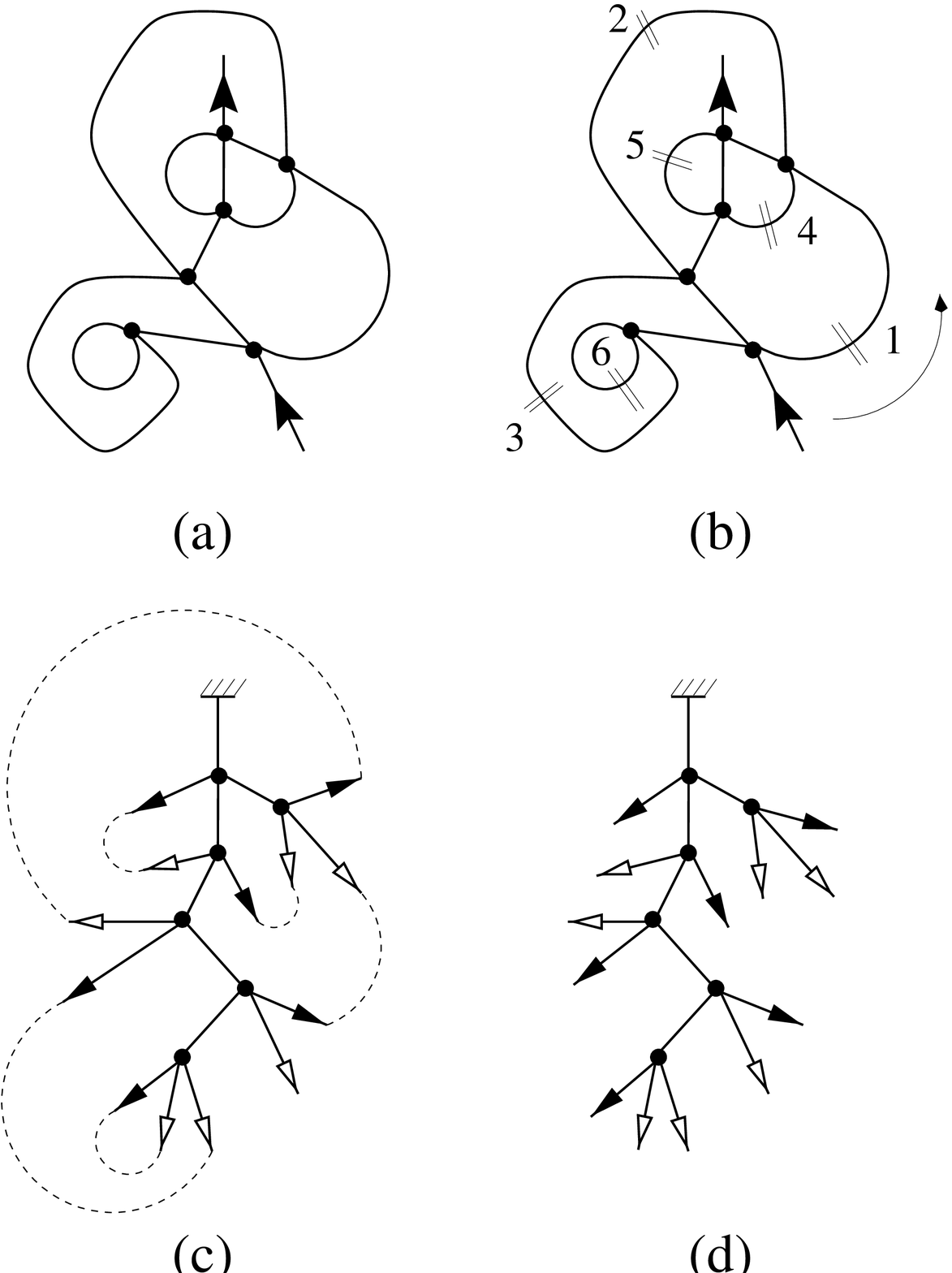}{10.cm}
\figlabel\bijecqua

We wish to compute the generating function $R\equiv R(g)$ for planar tetravalent
graphs with two distinguished (say in- and out-coming) univalent vertices, and with
a weight $g$ per tetravalent vertex. 
These will be referred to as ``two-leg diagrams" in the following, the legs
denoting simply the edges connecting the two univalent vertices to the graph.
For definiteness, we will always represent these planar graphs with the incoming
leg adjacent to the external face. Note that the outcoming one need not be
adjacent to the same face, as in the example of Fig.\bijecqua\ (a).

The computation of $R$ relies on the following bijection, illustrated in Fig.\bijecqua, 
between two-leg diagrams and so-called rooted blossom-trees \CENSUS.
Starting from a two-leg diagram (Fig.\bijecqua\ (a)), 
let us first visit all edges adjacent to the external face,
starting from the incoming leg and in counterclockwise order. Successively, each
visited edge is cut iff the cut diagram remains connected
(Fig.\bijecqua\ (b)). The cut ends  of the edge are then decorated
respectively with a black and a white leaf. 
Once all edges adjacent to the external face have been
visited, we repeat the procedure with the newly cut diagram. The process ends when all faces
of the original diagram have been merged with the external one. The resulting graph
is nothing but a planar tree (it has only one face), with black and white leaves (Fig.\bijecqua\
(c)).
Note that by construction
there is exactly one black leaf connected to each internal vertex of the tree. 
Finally, we replace the incoming vertex by a white leaf and the outcoming one by a root, so that there
is exactly one more white leaf than black ones (Fig.\bijecqua\ (d)).
We define rooted blossom-trees as planar rooted ternary trees with black and white leaves,
and such that there is exactly one black leaf at each internal vertex. Our cutting procedure
has produced a rooted blossom-tree out of any two-leg diagram. The process however is readily
seen to be invertible as there is a unique way of re-connecting the black-white leaf pairs
into edges, by connecting each black leaf to the first available white leaf in counterclockwise
order around the tree (the dashed lines of Fig.\bijecqua\ (c)). 
This establishes the desired bijection between the two objects.

Couting rooted blossom-trees is now an easy task, performed for instance by inspection
of all possible environments of the vertex connected to the root. This leads
straightforwardly to the relation
\eqn\ralaRqua{\eqalign{
&\figbox{8.cm}{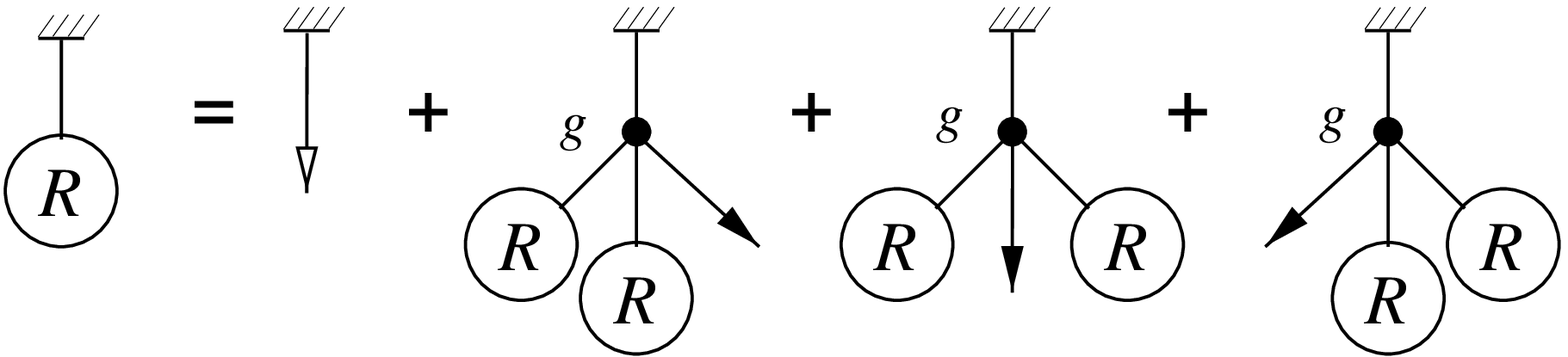} \cr
&\ \ \ \ \ R=1+3 g R^2 \cr} }
where the first term $1$ counts the possibility that the root be directly connected to a
leaf, and the second term accounts for the three possible positions for the black leaf around this vertex,
itself receiving the weight $g$, while the two other descendents of the vertex are themselves
rooted blossom-trees. From its very definition as counting function, $R$ admits a power series
expansion in $g$, with $R=1+O(g)$. This fixes it uniquely to be 
\eqn\Rfqua{ R={1-\sqrt{1-12g}\over 6g} }
The series for $R$ has a finite convergence radius $g_c=1/12$. When $g$ approaches $g_c$
(critical limit),
the contribution of large graphs becomes dominant, and we learn that the number of graphs
with $N$ vertices behaves as $g_c^{-N}/N^{3/2}$ for large $N$.

Note that $R(g)=C(3g)$ where $C$ denotes the generating function for Catalan numbers
$c_N={2N \choose N}/(N+1)$, which count among other things the rooted planar binary trees
with $N$ inner vertices, with the convention that $c_0=1$. The number
of rooted blossom trees with $N$ inner vertices
is obtained by considering rooted planar binary trees with all leaves white
and by decorating each vertex with a black leaf: it reads therefore $3^N c_N$ as there are 
three choices for the position of the black leaf at each inner vertex.

\subsec{General case of graphs with even valences}

The bijection of previous section may be extended to include two-leg-diagrams of graphs
with arbitrary even valences $v=4,6,8...$ We repeat the exact same cutting procedure
and end up with some generalized rooted blossom-trees, such that each inner vertex say of valence
$v=2k$ has exactly $k-1$ black leaves attached to it. The corresponding generating function
$R\equiv R(g_4,g_6,g_8,...)$ with say weights $g_{2k}$ per $2k$-valent vertex, $k=2,3,4...$
obeys the following relation
\eqn\relgeR{ R=1+\sum_{k\geq 2} g_{2k} {2k-1\choose k} R^k }
obtained again by inspecting all possible environments of the vertex attached to the root.
The combinatorial factor ${2k-1\choose k}$ accounts for the number of choices for the positions
of the $k-1$ black leaves around the vertex, while the remaining $k$ descendents are themselves
rooted blossom-trees.

Again, $R$ is the unique solution to \relgeR\ that admits a power series expansion
in the $g_{2i}$'s, with $R=1+O(g_{2i})$ for all $i\geq 2$.

\newsec{Planar graphs and trees II: arbitrary valences}

This section extends the results of the previous one to graphs with both even and odd valences.
The first consequence of allowing for odd valences is the existence of one-leg diagrams with
only one (outcoming) external leg. These will be represented in the plane like two-leg diagrams,
but the leg need not be adjacent to the external face.

\subsec{Trivalent graphs}

Let $S\equiv S(g)$ and $R\equiv R(g)$ denote the generating functions for respectively
one- and two-leg diagrams of trivalent planar graphs, with a weight $g$ per trivalent vertex.
Applying the cutting procedure of previous sections to one and two-leg diagrams, we end
up with two types of rooted blossom trees which we call S-trees and R-trees
respectively. Note that the unique (outcoming) leg of the one-leg diagrams is replaced by
a root in the corresponding blossom-tree. 
To characterize S- and R-trees, let us introduce the charge
$q$ of a tree as its number of white leaves minus that of black ones. For instance, the
blossom-trees of Sects. 2.1 and 2.2 above have all charge $q=1$, the same holds for the present R-trees,
while the S-trees are neutral, with charge $q=0$. Now  S-trees and R-trees are characterized among rooted
binary blossom-trees as having only descendent subtrees (not reduced to black leaves)
of charge 0 or 1, while their total
charge is $0$ and $1$ respectively (this is easily done by following the effect of the cutting procedure
on the original graph, see \CENSUS\ for details). The cutting procedure establishes a bijection between
one- and two-leg diagrams and S- and R-trees respectively.

The latter are easily counted, again by inspection of the local environment of the vertex attached
to the root. We get the coupled relations:
\eqn\coptri{\eqalign{&\figbox{6.cm}{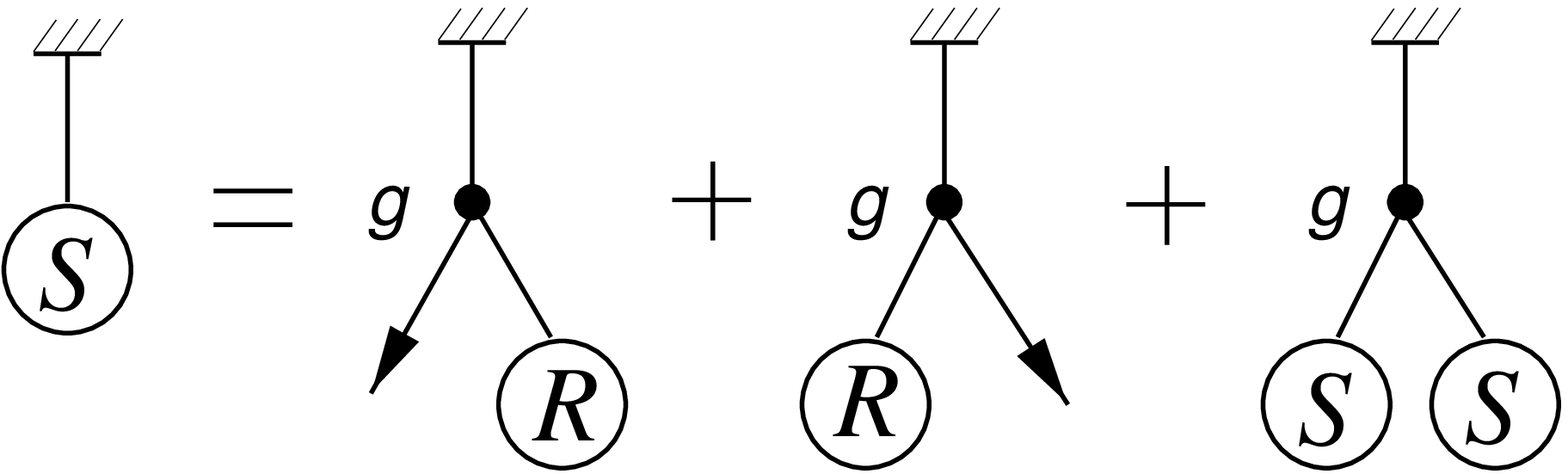}\cr
&\ \ \ \ \ S=2 g R + g S^2 \cr
&\figbox{6.cm}{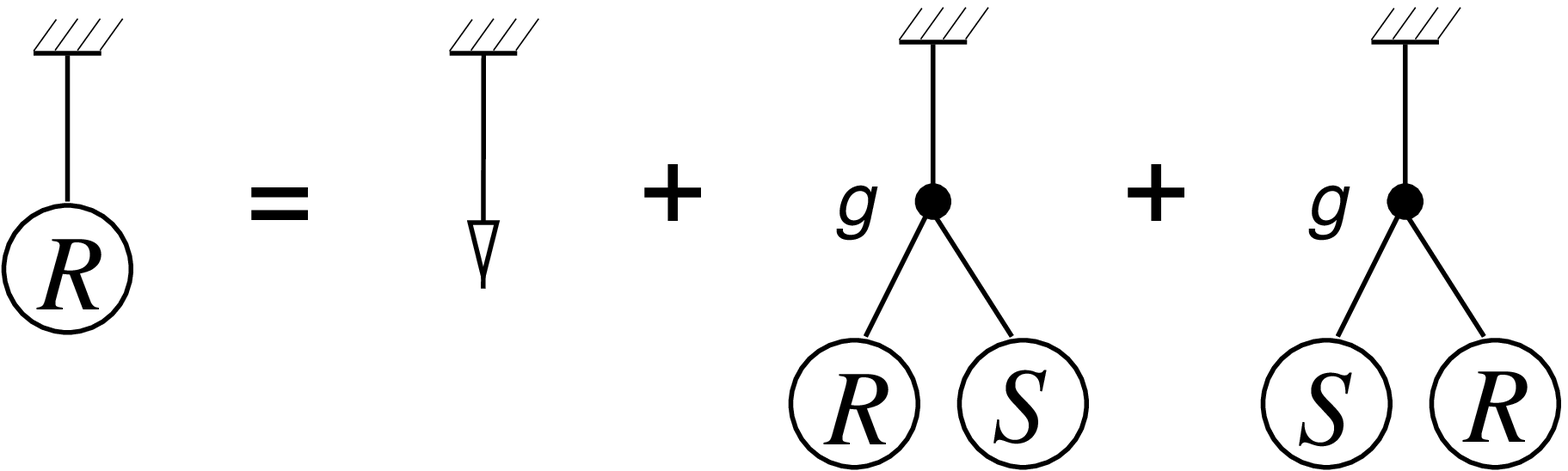}\cr
&\ \ \ \ \ R=1+2 g R S \cr}}
respectively displaying contributions from a vertex with total charge 0 (with one black leaf
with $q=-1$ and one descendent R-tree with $q=1$, and two possible positions for the black leaf,
or with two descendent S-trees), and from a vertex with total charge $1$ (with one descendent
S-tree with $q=0$ and one descendent R-tree with $q=1$, and two possible relative positions for
these). 
 
The generating functions $R,S$ are uniquely determined by the relations \coptri\ and the fact that
they admit power series expansions $R=1+O(g)$, $S=O(g)$.

\subsec{General case}

The trivalent case is easily extended to the case of arbitrary (even or odd) valences weighted
by $g_3,g_4,g_5,...$ per tri-, tetra-, penta-,... valent vertex, in which the very same
cutting procedure now leads to generalized rooted blossom S- and R-trees, now blossom trees
of arbitrary valences $v=3,4,5,...$ again further characterized by the fact that all
their descendent subtrees not reduced to a black leaf have charge $0$ or $1$, and by their total charge
0 and 1 respectively. This allows to count them straightforwardly, with the
coupled relations:
\eqn\coupodd{\eqalign{
S&=\sum_{k\geq 3} g_{k}  \sum_{j=0}^{[{k-1\over 2}]} {k-1\choose j}{k-1-j\choose j} R^j S^{k-1-2j} \cr
R&=1+\sum_{k\geq 3} g_{k} \sum_{j=0}^{[{k-2\over 2}]} {k-1\choose j}{k-1-j\choose j+1}
R^{j+1} S^{k-2-2j} \cr}}
where the combinatorial factors account for the possible ways of positioning $j$ black leaves,
$j$ or $j+1$ descendent R-subtrees and the remaining S-subtrees on the vertex attached to the root.
Again, eqs.\coupodd\ determine completely $R$ and $S$ with $R=1+O(g_{i})$ and $S=O(g_{i})$
for all $i\geq 3$.

For illustration, in the case of tri/tetravalent graphs, where only $g_3,g_4$ are non-zero, we
have the equations
\eqn\platritet{\eqalign{
S&=g_3(2R+S^2)+g_4(6RS+S^3)\cr
R&=1+2g_3RS+3g_4 R(R+S^2)\cr}}

\newsec{Planar graphs and trees III: bipartite graphs}

We now turn to the slightly more involved case of bipartite (i.e. vertex-bicolored, 
say black and white)
graphs. In the language of matrix models, these correspond to the case of two coupled
matrices.

\subsec{$p$-valent case}

Let us consider two-leg diagrams of vertex-bicolored $p$-valent planar graphs,
and their generating functions
with a weight $g$ (resp. $\tilde g$) per $p$-valent black (resp. white) vertex. 
We must also indicate the color of the 
vertices to which the in and out-coming legs are attached, 
and this leads to a priori distinct generating
functions. For simplicity, we restrict ourselves to only diagrams with incoming
(resp. outcoming) leg attached to a white (resp. black) vertex
with generating function $R\equiv \Rwb(g,{\tilde g})$, and also include the single graph
whose incoming and outcoming legs are directly attached to one-another, without vertex,
contributing $1$ to $R$. 
\fig{A sample bipartite $p$-valent two-leg diagram with $p=4$ (a) is iteratively cut
into a rooted bipartite blossom-tree (c) by the usual procedure, with the restriction
that only edges originating from a black vertex may be cut. The edges to be cut
are indicated in (b) with their order of visit, counterclockwise around the graph.}{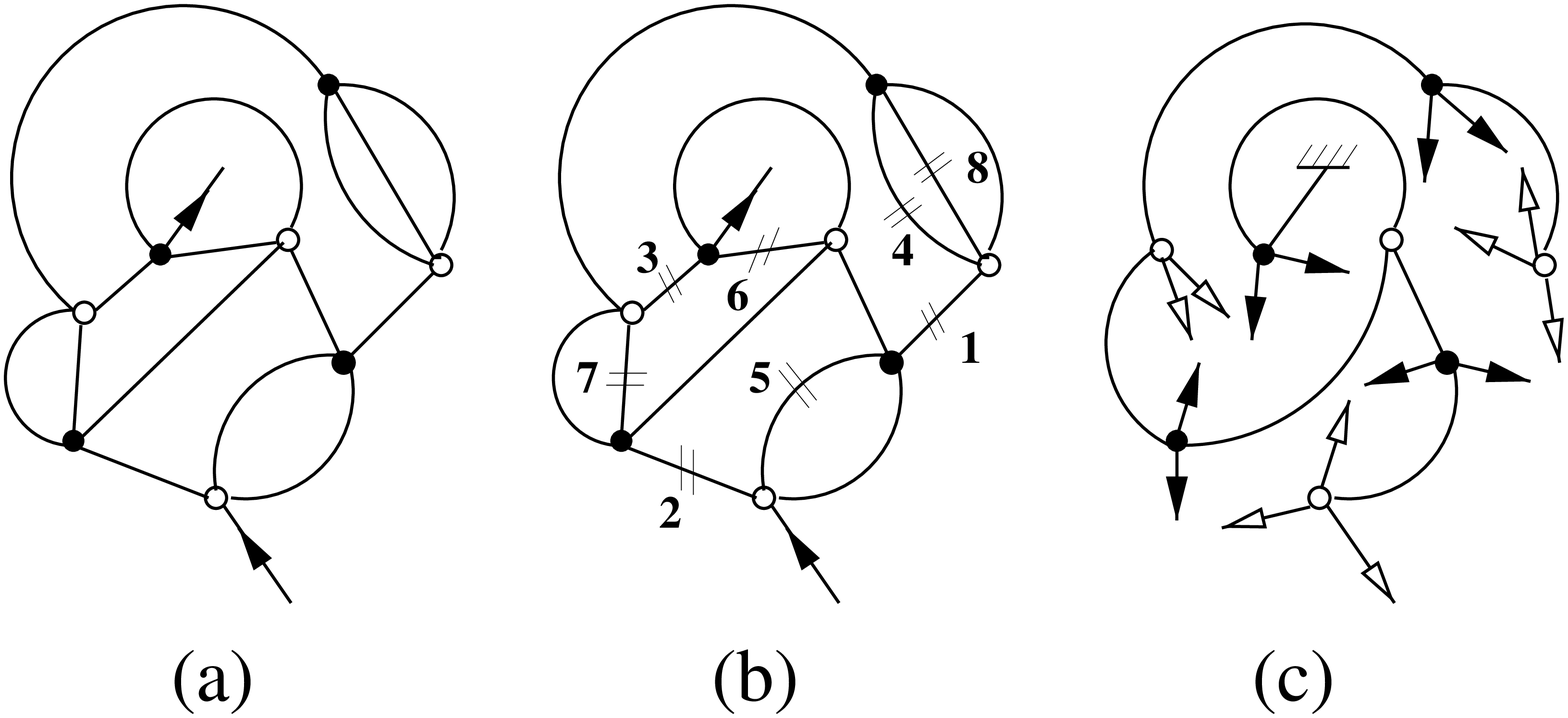}{12.cm}
\figlabel\pcolor
The above cutting procedure is now slightly modified, with the additional constraint that an edge 
may be cut only if it moreover 
originates from a black vertex (see Fig.\pcolor\ for an illustration in the case $p=4$). 
This leads to a new kind of blossom-trees (Fig.\pcolor\ (c)), with bicolored
vertices, and such that black or white leaves may only be connected to vertices of the same color,
while each black vertex has exactly $p-2$ black leaves attached to it. 
The tree still has the property of
having one more white leaf than black ones (i.e. a total charge
of $q=1$), as the incoming leg is replaced by a white leaf (which
is compatible with the above rule, as the incoming leg is attached to a $p$-valent white vertex). 
Moreover, the root (former outcoming vertex) is attached to a black vertex.
The generating function $R$ obeys the relations:
\eqn\obeytribi{ \eqalign{
R&=1+(p-1) g X \cr
X&={\tilde g}R^{p-1}\cr}}
The first line is obtained by inspection of all possible environments of the root: (i) it may simply have
one white leaf attached to it or (ii) it may have a black vertex attached to it, itself with 
$p-2$ black leaves and a descendent rooted blossom-tree of total charge $+(p-1)$ with a black root, 
and with generating function $X$.
The second line expresses these latter trees according to the environment of the 
white vertex attached to the root, having $p-1$ descendent rooted blossom-trees of charge $+1$, all 
generated by $R$.
We finally get
\eqn\pata{R=1+ (p-1)g {\tilde g} R^{p-1} }
Note that the generating function $C_p(x)$ for rooted $p$-valent planar trees with a weight $x$ per vertex
satisfies $C_p(x)=1+xC_p(x)^{p-1}$ and $C_p(x)=1+O(x)$. The corresponding number of trees with $N$ vertices
reads $C_N^{(p)}= {1\over 1+(p-1)N}{(p-1)N\choose N}$. These numbers are
also known as the Fuss-Catalan numbers, and reduce to the ordinary Catalan numbers for $p=3$. 
Finally the number of rooted blossom-trees with $N$ vertices is simply $R\vert_{g^N{\tilde g}^N}
=(p-1)^N C_N^{(p)}$. For large $N$, it behaves as $g_c^{-N}/N^{3/2}$, where
$g_c=(p-2)^{p-2}/(p-1)^p$.

\subsec{$p$-constellations}

\fig{The cutting procedure is applied to a two-leg diagram
of a 3-constellation (a), with incoming (resp. outcoming) leg attached to a white (resp. black)
vertex. The edges are visited in
counterclockwise order around the graph (b), and cut iff (i) this leaves the resulting graph connected
(ii) the cut edge originates from a black vertex. We have indicated the chronological order of the
cut edges, from 1 to 6. We finally replace each cut edge by a pair of black/white leaves,
and the incoming (resp. outcoming) leg by a white leaf (resp. root).}{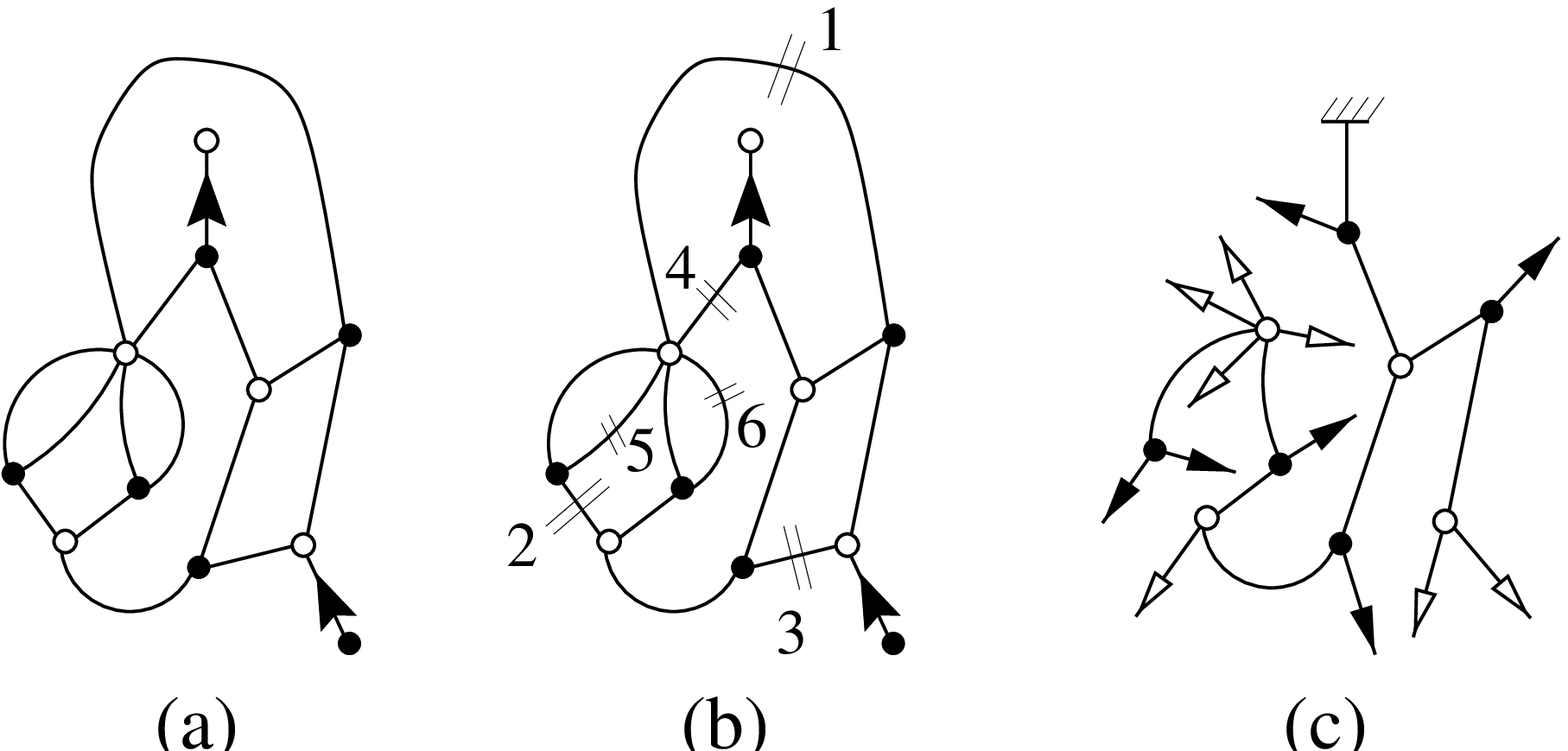}{12.cm}
\figlabel\constellation

The
$p$-constellations, introduced in \BMS, 
are vertex-bicolored (black and white) graphs such that say black vertices have 
all valence $p$, while white ones may have valences arbitrary multiples of $p$.  
We again consider two-leg diagrams of such planar graphs with say incoming leg connected to a 
white vertex and outcoming leg connected to a black vertex (see Fig.\constellation\ (a)
for an illustration with $p=3$), or both legs being directly 
connected without vertex,
and count them with a weight $g$ per ($p$-valent) black vertex and weights ${\tilde g}_m$
for $mp$-valent white vertices, $m=1,2,3...$ Let $R=\Rwb(g;\{{\tilde g}_i\})$ denote the
corresponding generating function.
The cutting procedure, illustrated in Fig.\constellation, remains the same as in the 
previous section, and leaves us with rooted
blossom-trees with bicolored vertices of total charge $+1$, such that leaves may only be 
connected to vertices of their own color, and that descendent subtrees may be of two types:
if their first vertex is black, they have total charge $1$; if its is white, they have total
charge $p-1$, the descendents of the root vertex of the latter being themselves 
only blossom trees of charge
$1$ or bunches of $p-1$ black leaves attached to a black vertex. 

By inspection of all possible local environments of the vertex attached to the root,
we may derive the following relations for $R$:
\eqn\saticonst{\eqalign{
R&=1+(p-1)g X\cr 
X&=\sum_{m\geq 1}{\tilde g}_m {mp-1\choose m-1}Y^{m-1} R^{(p-1)m} \cr
Y&=g}} 
where $R$, $X$, $Y$ generate rooted blossom-trees of respective charges $1$, $p-1$, $1-p$,
starting respectively with a black, white, black vertex. $Y=g$ is due to the fact that
the only tree contributing to $Y$ has a black vertex and $p-1$ black leaves attached to it.
In the case $p=3$, this reads
\eqn\troiconst{\eqalign{
&\figbox{5.cm}{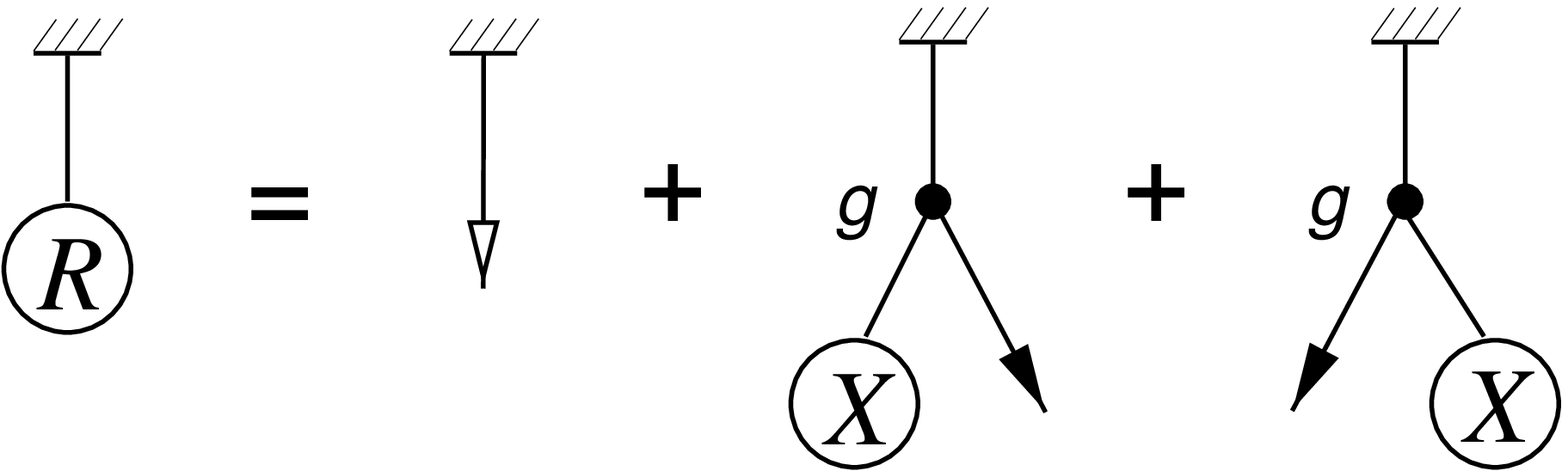}\cr
&\ \ \ \ \ R=1+2g X\cr
&\figbox{13.cm}{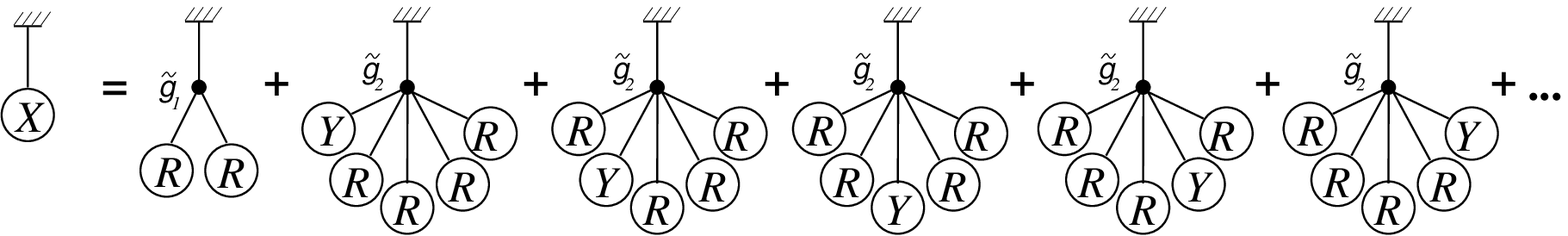} \cr
&\ \ \ \ \ X={\tilde g}_1 R^2 +5 {\tilde g}_2 Y R^4+ \cdots  \cr
&\figbox{2.cm}{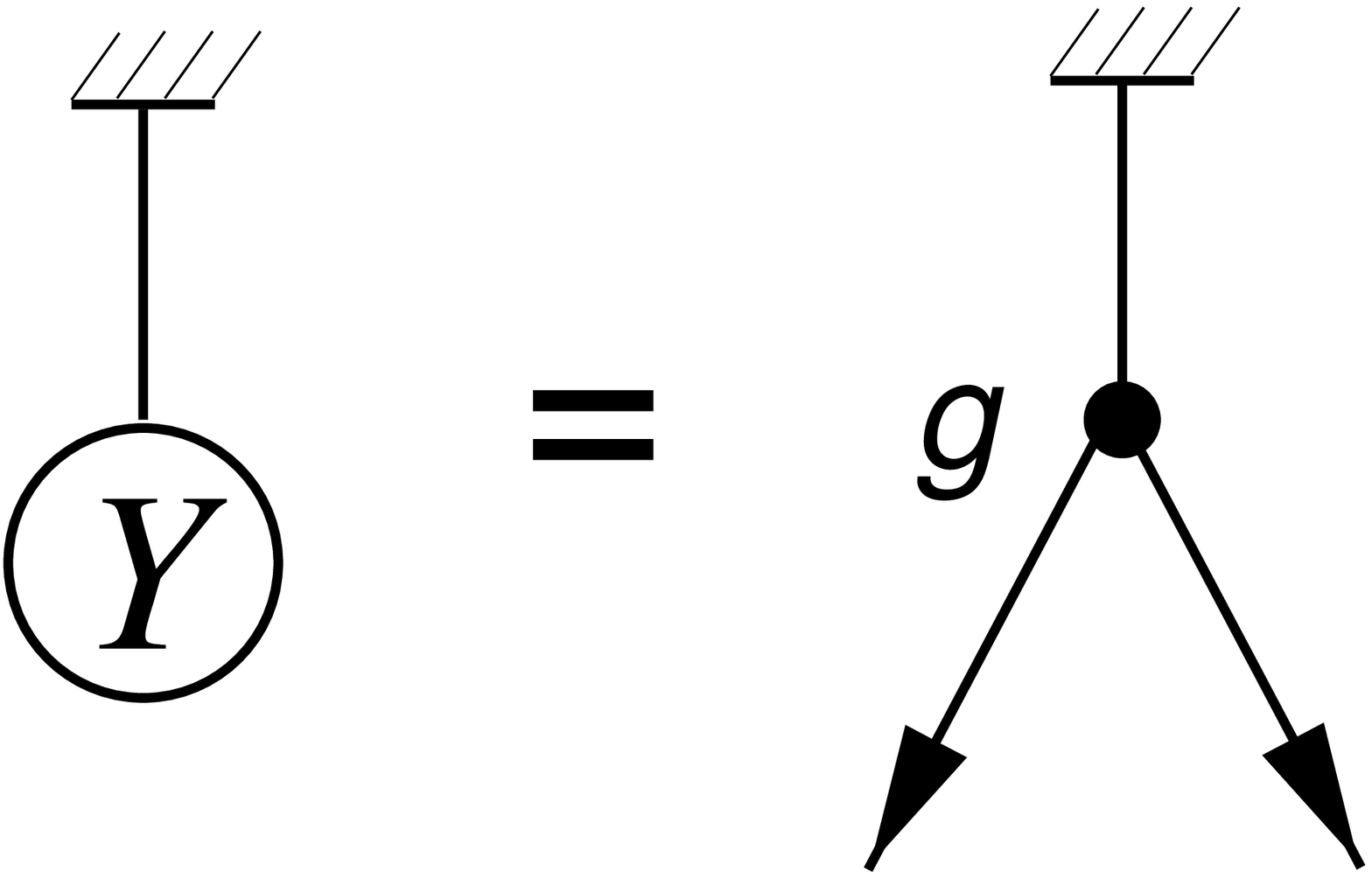} \cr
&\ \ \ \ \ Y=g\cr}}

Eliminating $Y$ and $X$ from \saticonst, we arrive at the single algebraic equation
\eqn\singal{ R=1+(p-1)\sum_{m\geq 1}{\tilde g}_m {mp-1\choose m-1}g^{m-1} R^{(p-1)m} }
$R$ is the unique solution to this equation such that $R=1+O(g,{\tilde g}_i)$.

Note that in the case $p=2$ of 2-constellations, we recover the even-valent graph result
\relgeR, upon taking $g_{2k}={\tilde g}_k$ for $k\geq 2$, while ${\tilde g}_1=0$ and $g=1$.
Indeed, in 2-constellations, the (2-valent) black vertices may be viewed as decorations of the edges
of an arbitrary graph with only white even-valent vertices. Setting $g=1$ precisely allows
to forget about these decorations, while ${\tilde g}_1=0$ simply eliminates 2-valent white vertices.

\subsec{Planar bipartite graphs and the Ising model}

Constellations are easily tractable objects, essentially due to the
triviality of one type of vertex (black here), whose valence remains
fixed. More generally, we would like to consider in all generality
vertex-bicolored graphs of arbitrary even valences for vertices of both colors.
For the sake of simplicity, we will restrict ourselves to two-leg diagrams in which
both legs are attached to white vertices.
Let us denote by $R\equiv \Rww(\{g_{2i}\};\{{\tilde g}_{2i}\})$ the generating function for
two-leg diagrams of planar bipartite graphs with weights $g_{2i}$ (resp. ${\tilde g}_{2i}$)
per $2i$-valent black (resp. white) vertex, $i=1,2,...$, and such that both
legs are attached to white vertices. 

Applying the now usual cutting procedure to two-leg diagrams, we end up with blossom-trees
of total charge $+1$ with bicolored vertices of arbitrary even valences. 
Their characterization however
is quite delicate, as it involves describing all their possible descendent subtrees.
These come in two forms according to the color of their vertex attached to the root:
if the latter is white, the possible descendent subtrees are either reduced to a black leaf
(of charge $-1$) or rooted vertex-bicolored
blossom trees of total charges $1$, $3$, $5$, ... ; if it is black, the possible descendent
subtrees are rooted vertex-bicolored
blossom trees of charges $1$, $-1$, $-3$, $-5$, ...   
Let us introduce the generating functions $R_i$ for vertex-bicolored
blossom-trees of total charge $2i-1$, $i=1,2,3,...$ and whose root is attached to a white vertex, 
together with the generating functions $X_i$ for vertex-bicolored
blossom-trees of total charge $1-2i$, $i=1,2,3,...$ and whose root is attached to a black vertex,
while $V\equiv X_{0}$ generates  vertex-bicolored
blossom-trees of total charge $1$ and whose root is attached to a black vertex, or the tree made of a single leaf
attached to the root, without any vertex (contributing $1$ to $V$). We also introduce the generating
function $R_{0}=1$ for the tree made of a single black leaf attached to the root, without any vertex.  
We now simply have to enumerate all possible environments of the vertex attached to the
root of each of these trees, according to the type of its attached
descendent subtrees.
This gives the following system
\eqn\folsis{\eqalign{
V&=1+\sum_{k\geq 1} g_{2k} \sum_{j_1,j_2,...,j_{2k-1}\geq 0\atop \Sigma j_l=k } 
R_{j_1}R_{j_2}...R_{j_{2k-1}}\cr 
X_m&=\sum_{k\geq 1} g_{2k} \sum_{j_1,j_2,...,j_{2k-1}\geq 0\atop \Sigma j_l=k-m} 
R_{j_1}R_{j_2}...R_{j_{2k-1}},\ \  m=1,2,3,...\cr
R_m&=\sum_{k\geq 1} {\tilde g}_{2k} \sum_{j_1,j_2,...,j_{2k-1}\geq 0\atop \Sigma j_l=k-m}
X_{j_1}X_{j_2}...X_{j_{2k-1}}, \ \ m=1,2,3,...\cr}}  
The desired generating function $R=R_1$ is the unique solution to this system where all $R$'s, $X$'s and $V$
admit power series expansions of the $g$'s and ${\tilde g}$'s.

The case of the Ising model on planar tetravalent graphs may be viewed as a particular 
case of the above, in which only $g_2,g_4,{\tilde g}_2,{\tilde g}_4$ are non-zero. 
To see this, recall that the Ising model on tetravalent planar graphs is defined by 
say coloring the vertices of an arbitrary tetravalent planar graph in black or white
(colors stand here for the spin up or down),
and counting the configurations with different ``nearest neighbor interaction" weights for edges connected
to vertices of the same color (weight $e^K$) or of different colors (weight $1$), while
black (resp. white) vertices are counted with a weight $g e^{H}$ (resp. $g e^{-H}$). Here $K$ and
$H$ are respectively the spin coupling and the external magnetic field of the Ising model. 
To make the contact with our model, we just have to resum all possible configurations obtained by adding 
arbitrary numbers of 2-valent black and white vertices on the edges of Ising configurations, in such a way
that bicoloration is restored. This entails adding any chain of black, white, black, ..., white
2-valent vertices between any white and black tetravalent vertices connected by an edge, or
any chain of alternating white, black, ..., white 2-valent vertices between any two black tetravalent 
vertices connected by an edge  or else any chain of black, white, ..., black 2-valent vertices
between any two white tetravalent vertices connected by an edge. Doing the resummations within
the configurations of our bicolored graphs produces an effective edge interaction
weight $w_{ab}$ according to the colors $a,b$ of the adjacent vertices: 
\eqn\interw{\eqalign{ w_{\n \n}&= {{\tilde g}_2\over 1-g_2{\tilde g}_2} 
=\figbox{5.cm}{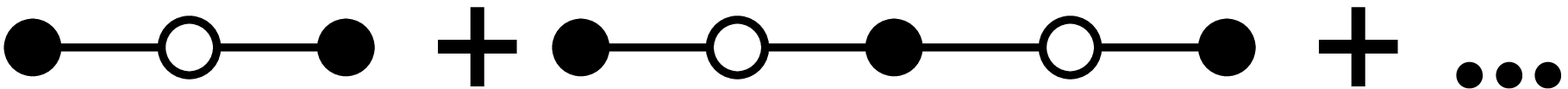} \cr 
w_{\b \b}&= {g_2\over 1-g_2{\tilde g}_2} 
=\figbox{5.cm}{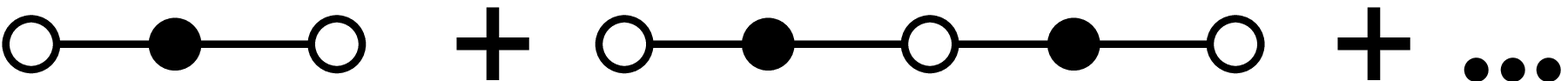} \cr
w_{\n \b}&=w_{\b \n}= {1\over 1-g_2{\tilde g}_2} 
=\figbox{4.5cm}{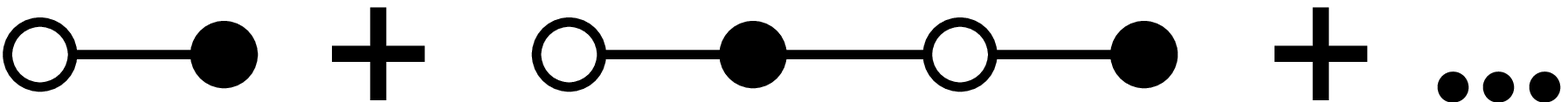} \cr}}
hence to identify our model with the Ising one,
we must take $g_2={\tilde g}_2= e^K$, while $g_4=(1-e^{2K})^2ge^H$ and
${\tilde g}_4=(1-e^{2K})^2ge^{-H}$, and the external legs must receive the extra weights 
$1/\sqrt{1-e^{2K}}$ each.

Restricting to the
symmetric case $g_2={\tilde g}_2\equiv c$
and $g_4={\tilde g}_4\equiv g$, the above equations simply read
\eqn\isiread{\eqalign{
&\figbox{13.cm}{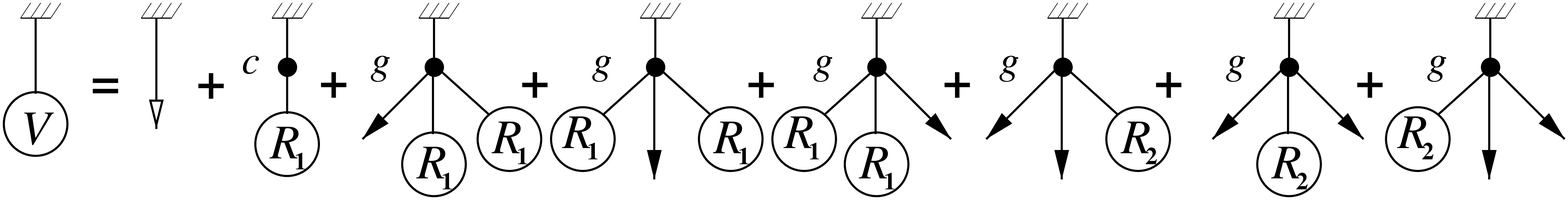}\cr
&\ \ \ \ \ V=1+c R_1+3 g R_1^2+3 g R_2 \cr
&\figbox{7.cm}{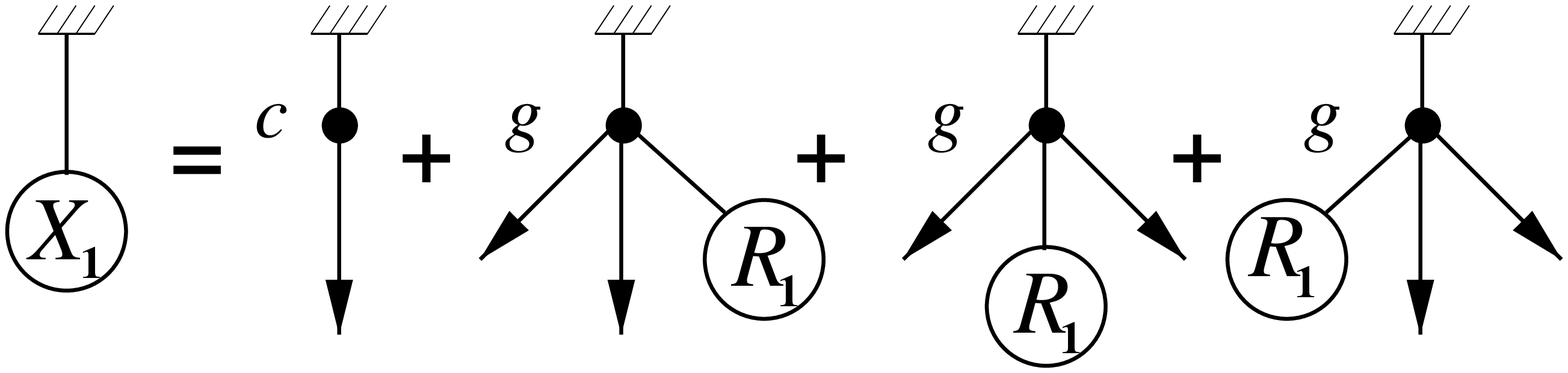}\cr
&\ \ \ \ \ X_1= c +3g R_1\cr
&\figbox{3.cm}{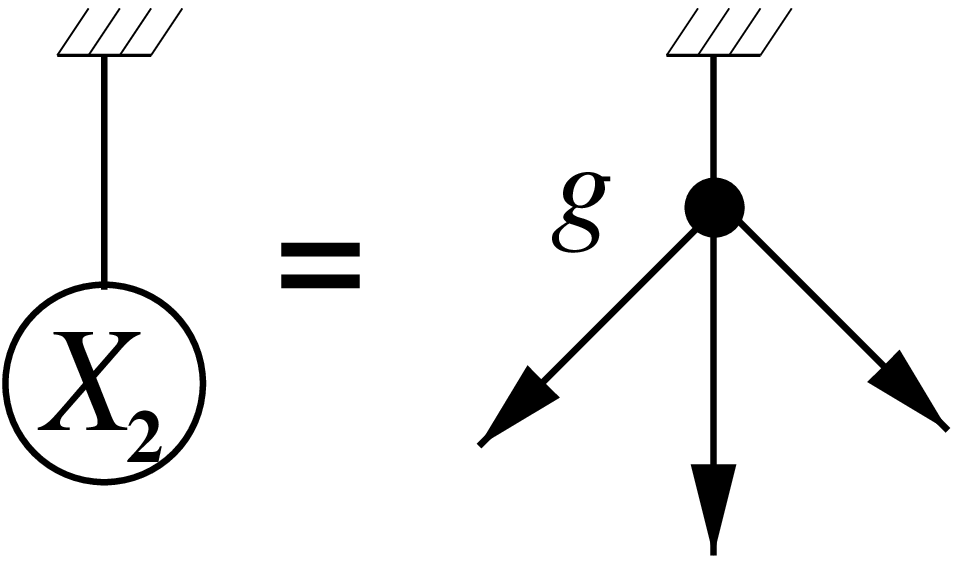}\cr
&\ \ \ \ \ X_2= g\cr
&\figbox{9.cm}{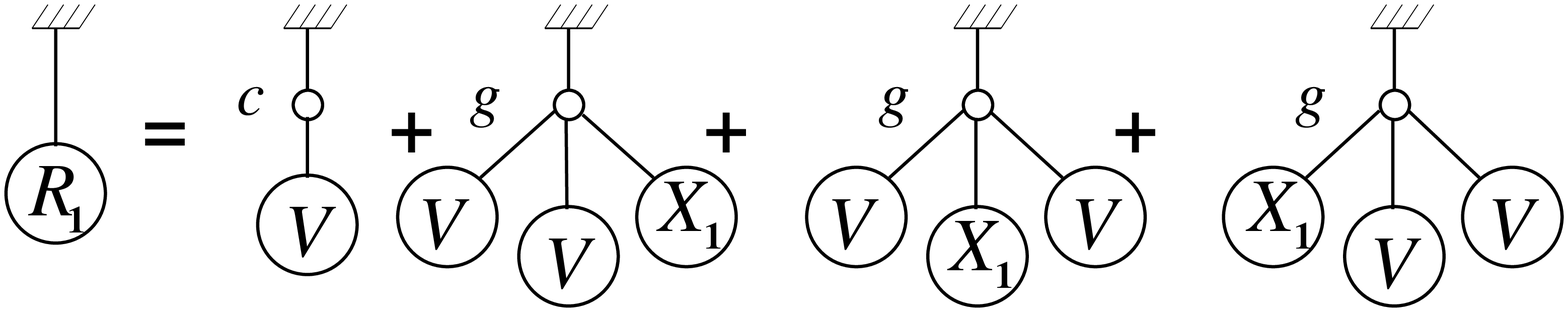}\cr
&\ \ \ \ \ R_1= c V+ 3g V^2 X_1\cr
&\figbox{3.cm}{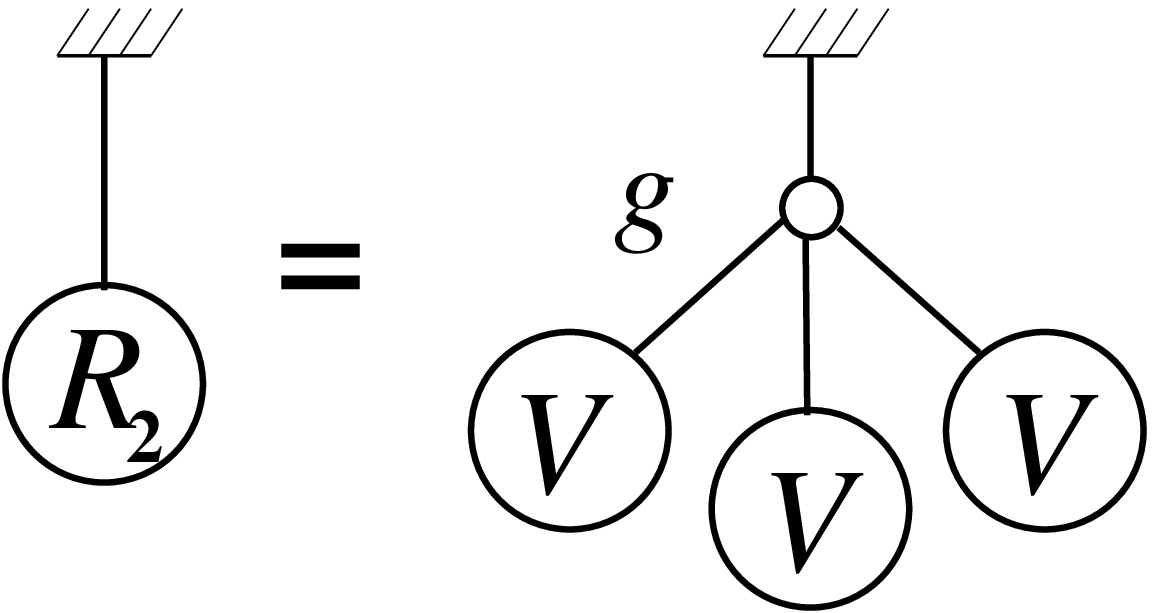}\cr
&\ \ \ \ \ R_2= g V^3 \cr}}
Eliminating $R_2$ and $X_1$, we get 
\eqn\finroneis{ V={R\over c+3 g R} \qquad {\rm with}\qquad 
R(c+3gR)^2(1-(c+3gR)^2)=(c+3gR)^3+3g^2R^3}
Here, $R\equiv R_1$ is the generating function for two-leg diagrams of planar tetravalent graphs
with Ising (black or white) spins decorating their vertices, with interaction weights
$w_{\n \n}=w_{\b \b}=c$, $w_{\n \b}=w_{\b \n}=1$ and a weight $g/(1-c^2)^2$ per tetravalent vertex,
and such that the two legs are attached to univalent black vertices, themselves weighted by $1/\sqrt{1-c^2}$.

\newsec{Geodesic distances}

In the previous sections, we have enumerated various one- and two-leg diagrams by establishing
bijections with suitable classes of blossom-trees. Note that in the plane representation we have chosen,
the face $F_1$ adjacent to the unique leg for one-leg diagrams is not necessarily the external face
$F_0$.
Accordingly in the case of two-leg diagrams, the face $F_1$ adjacent to the out-coming leg need not be the 
external face $F_0$, itself adjacent to the in-coming leg. 
In this section, we show how to keep track of the {\it geodesic distance} 
between $F_0$ and $F_1$, namely the smallest number of edges to be crossed in a path from $F_0$ to $F_1$.
By a slight abuse of language, this distance will also be referred to as the distance between the
legs.

\subsec{Keeping track of the geodesic distances}

The main feature of the previous sections is a sort of unified formulation of planar graphs 
in the language of blossom-trees. Note that going back from blossom-trees to graphs is a straightforward
step, as there is a unique way of reconnecting the black and white leaves into edges: this is done 
by simply connecting each black leaf to the first available white leaf in counterclockwise direction.
This process leaves us in the case of trees of charge $1$ with exactly one unmatched white leaf,
which is taken as outcoming leg, while the root is the incoming one. For trees of charge 0, all pairs
are exhausted and only the root remains as uniqe leg. In both cases, we note that the geodesic distance between
the faces $F_0$ (external) and $F_1$ (adjacent to the former root) is simply given by the number of
black-white edge pairs that separate the root from the external face after recombination. 
Keeping track of this geodesic distance simply amounts to keeping track of the black leaves ``in excess"
that require encompassing the root to be connected to their white {\it alter ego}.

\fig{The contour walk of a rooted blossom-tree. Visiting the tree (a) in 
clockwise direction starting from the root,
one keeps a record of the type of leaves encountered in the form of a walk on the integer line (b), 
starting at the origin,
and with steps up (for a black leaf)
and down (for a white one). The maximum reached by the walk is nothing but the geodesic distance
separating the root from the external face in the recombined graph. This distance is 3 in the
present case, as one readily checks by closing the tree back into a planar graph
(with tetra/hexavalent vertices here). Concentrating on the environment of the vertex
attached to the root, we see that each descendent subtree corresponds to a portion of
the walk (c), with a certain relative maximum. Expressing the global maximum of the
contour walk in terms of the relative maxima of its portions allows for writing
a recursion relation for the generating function for rooted blossom-trees whose
root is at a maximum distance $n$ from the external face of the recombined graph.}{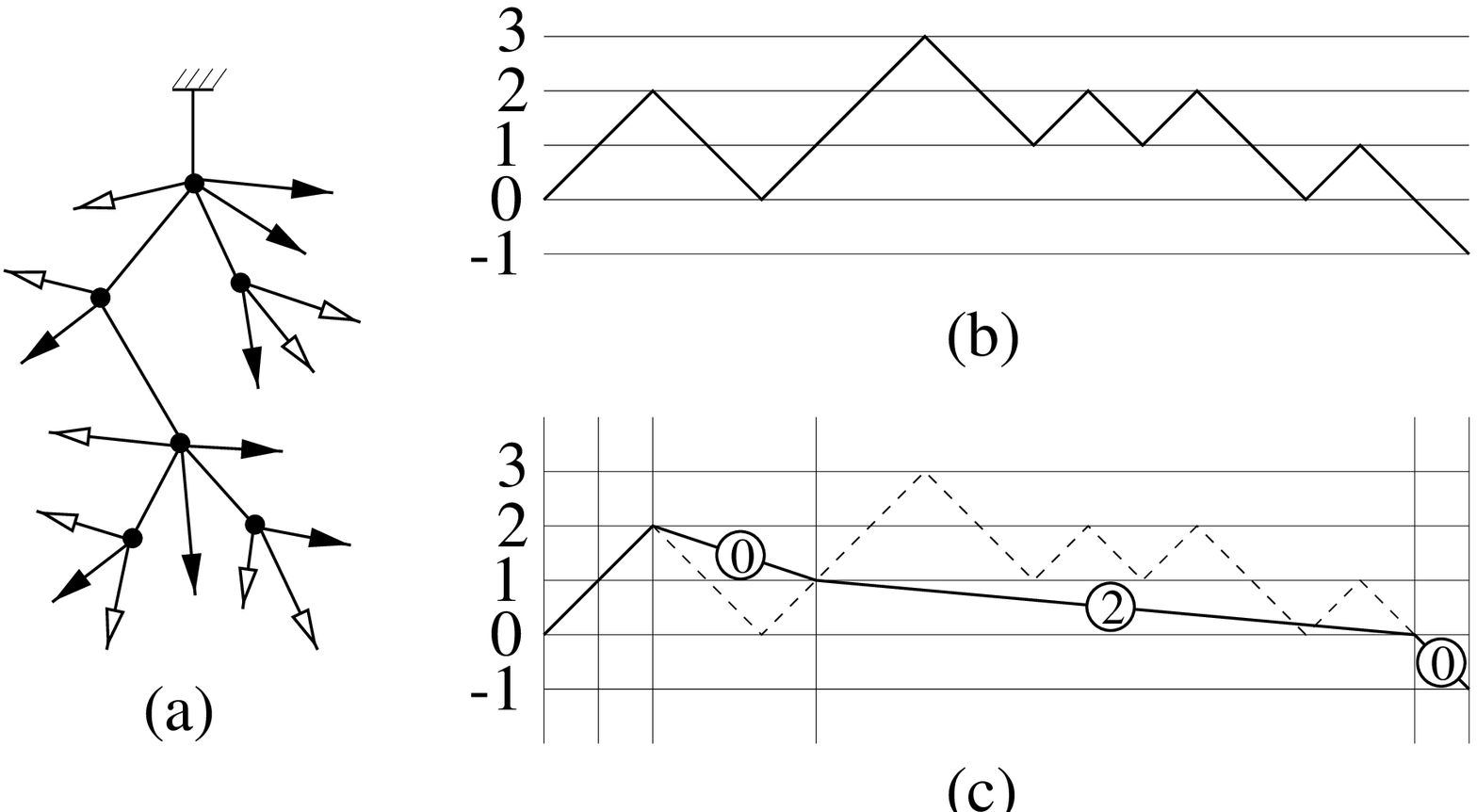}{13.cm}
\figlabel\contour

This is done in all generality by attaching to each blossom-tree a ``contour walk", namely a walk on
the relative integer line with steps $\pm 1$, obtained as follows (see Fig.\contour\
for an example). One starts from the root
of the blossom-tree and visits in clockwise direction all leaves around the tree. Starting from the
coordinate $0$, we make a step $+1$ (resp. $-1$) for each encountered black (resp. white) leaf. 
In a blossom-tree of charge $k$, such a walk will end up at coordinate $-k$. Now the number of excess black
leaves responsible for the geodesic distance between $F_1$ and $F_0$ is simply the maximum coordinate
reached by the contour walk.

In view of this result,  and of the form of all relations determining the blossom-trees of the previous sections,
it is natural to introduce by analogy with the generating functions say $X$ for some particular type of blossom-trees
the generating function $X_n$ for the same blossom-trees with a geodesic distance of at most $n$
between $F_1$ and $F_0$. To obtain from $X_n$ the generating function for blossom-trees with geodesic distance
equal to $n$ betwen $F_0$ and $F_1$, we simply have to take the difference $X_n-X_{n-1}$.
Keeping track of $n$ then boils down to expressing the maximum of the contour walk of a tree in terms of
those of the individual contour walks of the blossom-trees descending from the vertex attached to the root,
following the same inspection procedure as before.  
This is done case by case in the following sections.

\subsec{Even valences}

Let us start with the tetravalent case. We find that eq. \ralaRqua\ must be transformed into
\eqn\quaca{ R_n=1+gR_n(R_{n-1}+R_n+R_{n+1}) }
where $R_n$ denotes the generating function for tetravalent planar graphs with two legs at
distance at most $n$, and weight $g$ per vertex.
The three terms on the r.h.s. correspond to respectively the black leaf on the right, in the middle or 
on the left
of the two other descendents of the vertex attached to the root (see the picture of eq.\ralaRqua). 
It is clear that the presence of the black leaf acts as a shift by $-1$ on the local distance $n$
while $R_n$ accompanies a shift by $+1$, when going clockwise around the vertex.

To have a compact notation for the general result, let us introduce a formal orthonormal basis
$|n\rangle$, $n\in \IZ$, with $\langle m|n\rangle =\delta_{m,n}$,
and an operator $\sigma$ acting as a shift $\sigma |n\rangle =|n+1\rangle$,
and its formal inverse $\sigma^{-1}$ such that $\sigma^{-1}|n\rangle=|n-1\rangle$. We introduce the operator
\eqn\opQ{Q=\sigma +\sigma^{-1}{\hat r}}
where $\hat r$ simply acts diagonally as ${\hat r}|n\rangle=R_n |n\rangle$. 
Note that the shift $\sigma$ may be represented by a black leaf, 
while its inverse always accompanies an $R$.
Then eq. \quaca\ takes the form 
\eqn\newquaca{ 1=\langle n-1 | (Q-gQ^3) |n \rangle }
More generally, in the case of arbitrary even valences, we have to write
\eqn\arbeven{  1=\langle n-1 | (Q-\sum_{k\geq 2} g_{2k} Q^{2k-1}) |n \rangle }
For instance in the case of tetra- and hexa-valent graphs, eq. \newquaca\ reads explicitly
\eqn\hexaqua{\eqalign{ R_n&=1+g_4 R_n(R_{n+1}+R_n+R_{n-1})+g_6 R_n\big(R_{n+1}R_{n+2}+R_{n+1}R_{n-1}\cr
&+R_{n-1}R_{n-2}+
R_{n+1}^2+R_{n-1}^2+R_n(2R_{n+1}+R_n+2R_{n-1})\big)\cr} }

A remark is in order. The equations \quaca-\hexaqua\ are valid only for $n\geq 0$, provided we use
$R_{-k}=0$, $k=1,2,...$ wherever they occur in the r.h.s. With these boundary conditions, and the
general fact that $R_n$ possess a power series expansion in $g$, with $R_n=1+O(g)$ for $n\geq 0$,
all $R_n$'s are then uniquely determined by \arbeven\ order by order in $g$. 

\subsec{Arbitrary valences}

Let $S_n$ (resp. $R_n$) denote the generating function for one- (resp. two-) leg diagrams of planar graphs
with arbitrary valences with weights $g_i$ per $i$-valent vertex, 
and such that $F_0$ (the external face) and $F_1$ (the face adjacent to the unique (resp. outcoming) leg) 
are distant by at most $n$.
For trivalent graphs, we find that eqs. \coptri\ must be transformed into
\eqn\transfortri{ \eqalign{
S_n&=g(R_n+R_{n-1})+g S_n^2 \cr
R_n&=1+g R_n(S_{n+1}+S_n)\cr}}
where in addition to the situation of previous section we simply note that $S$'s don't affect the distance
counting (no shift). This suggests to introduce in the general case of arbitrary valences the operator
\eqn\moregenQ{Q=\sigma+\sigma^{-1}{\hat s} \sigma+\sigma^{-1}{\hat r}}
where $\hat s$ acts diagonally on the basis $|n\rangle$
as ${\hat s}|n\rangle =S_n |n\rangle$, and 
in terms of which we simply have to write
\eqn\intow{\eqalign{
0&= \langle n| (Q-\sum_{i\geq 3} g_i Q^{i-1}) |n \rangle \cr
1&= \langle n-1| (Q-\sum_{i\geq 3} g_i Q^{i-1}) |n \rangle \cr}}
For tri- and tetra-valent graphs, this reads
\eqn\tritetr{\eqalign{
S_n&= g_3(R_n+R_{n-1}+S_n^2) +g_4 (R_{n}(S_{n+1}+2S_n)+R_{n-1}(S_{n-1}+2S_n)+S_n^3) \cr
R_n&=1+ g_3R_n(S_n+S_{n+1}) + g_4 R_n(S_n^2+S_nS_{n+1}+S_{n+1}^2+R_{n+1}+R_n+R_{n-1}) \cr}} 

Note that when $g_3=0$, we find the solution $S_n=0$ identically, and eq. \tritetr\ reduces to
the tetravalent case \quaca. More generally, imposing that all odd $g$'s vanish leads to the
solution $S_n=0$ (as there are no one-leg diagrams with only even valences), 
and we recover the even-valent case of previous section. 

\subsec{Constellations}

Let $R_n$ denote the generating function for two-leg diagrams
of $p$-constellations with a weight $g$ per white 
($p$-valent) vertex, and weights ${\tilde g}_i$ per black $pi$-valent vertex, $i=1,2,...$,
whose incoming (resp. outcoming) leg is attached to a black (resp. white) vertex 
(or both are connected without vertex and the corresponding unique graph contributes $1$ to $R_n$ for all
$n$) and such that the geodesic distance between the two legs is at most $n$.
Here the notion of geodesic distance is defined according to the
rules used in the cutting procedure of Sect.5.3, namely the geodesic distance between $F_0$ and $F_1$ 
is the minimal number of edges to be crossed in a path going from $F_0$ to $F_1$, 
and such that at each edge-crossing the white vertex is always on the right. 
Following the same reasoning as in the previous sections, we are now led to 
the introduction of two operators $Q_1$ and $Q_2$ which generate, upon taking powers,
the successive decorations of the vertex attached to the root, respectively in the case
of a black and white vertex. We have
\eqn\qops{\eqalign{Q_1&=\sigma+\sigma^{-1} {\hat x} \sigma^{2-p} \cr
Q_2&=\sigma^{-1}{\hat r}+\sigma^{p-2} {\hat y} \sigma \cr}}
where the shift $\sigma$ represents a
single black leaf, while $\hat x$ represents rooted blossom-trees of charge $p-1$ whose first vertex is
white. Again, ${\hat x},{\hat y},{\hat r}$ act diagonally on the basis $|n\rangle$
with eigenvalues $X_n,Y_n,R_n$ respectively.
The equations \saticonst\ now become 
\eqn\saticonstn{\eqalign{
1&=\langle n-1|(Q_2-gQ_1^{p-1})|n\rangle \cr
0&=\langle n+p-1|(Q_2-gQ_1^{p-1})|n\rangle \cr
0&=\langle n-p+1|(Q_1-\sum_{i\geq 1} {\tilde g}_i Q_2^{pi-1})|n\rangle\cr}}
In the particular case of 3-constellations, with say only $g,{\tilde g}_1,{\tilde g}_2$ non-zero,
these read for instance
\eqn\instcons{\eqalign{
R_n&=1+g (X_n+X_{n-1})\cr
Y_n&=g\cr
X_n&={\tilde g}_1 R_n R_{n+1} +{\tilde g}_2 g R_nR_{n+1}(R_{n+3}R_{n+2}+R_{n+2}R_{n+1}+R_{n+1}R_n
+R_{n}R_{n-1}) \cr
}}
while in the case of only (black and white) $p$-valent vertices of Sect. 4.1, eqs.\saticonstn\
read
\eqn\onlybw{\eqalign{
R_n&=1+g (X_n+X_{n-1}+...+X_{n-p+2})\cr 
X_n&={\tilde g}_1 R_n R_{n+1}...R_{n+p-2} \cr
}} 
for the generating function $R_n$ for two-leg-diagrams of bipartite $p$-valent 
graphs with incoming leg attached to a white vertex and outcoming leg attached
to a black one, and such that the geodesic distance from the in- to the out-coming leg
is at most $n$.

\subsec{Bipartite even-valent graphs and the Ising model}

In the general case of bipartite even-valent graphs, we are led to the introduction of two
operators $Q_1,Q_2$ with the following structure
\eqn\strucqonetwo{\eqalign{
Q_1&= \sigma +\sum_{k\geq 1} \sigma^{1-2k} {\hat r}^{(k)} \cr
Q_2&= \sigma^{-1} {\hat v} +\sum_{k\geq 1} \sigma^{-1} {\hat x}^{(k)} \sigma^{2k} \cr}}
where the ${\hat r}^{(k)}, {\hat x}^{(k)}, {\hat v}$ all act diagonally 
with eigenvalues $R_n^{(k)}, X_n^{(k)}, V_n$.
The latter are nothing but the generating functions for sets of 
rooted blossom-trees restricted by $n$, respectively
starting with a white, black, black vertex, and with charges $2k-1,1-2k,1$ respectively.
We are actually interested in computing $R_n\equiv R_n^{(1)}$, the
generating function for two-leg diagrams of bipartite planar graphs with weights
$g_{2i},{\tilde g}_{2i}$ per $2i$-valent black, white vertex, such that moreover the two legs are
distant by at most $n$. Again, the distance from $F_0$ to $F_1$ is defined as the minimal number 
of edges to be crossed in a path from $F_0$ to $F_1$, such that at each edge crossing 
the white vertex is always on the right.
This definition allows to keep track of this distance on the trees themselves, as the number
of ``excess" black leaves which upon recombination with white ones encompass the root of the tree.
The equations determining $R_n$ are simply
\eqn\eqdert{\eqalign{
1&= \langle n-1| (Q_2 -\sum_{i\geq 1} g_{2i} Q_1^{2i-1}) | n\rangle \cr
0&= \langle n+2m-1| (Q_2 -\sum_{i\geq 1} g_{2i} Q_1^{2i-1}) | n\rangle, \qquad m=1,2,... \cr  
0&=  \langle n-2m+1| (Q_1 -\sum_{i\geq 1} {\tilde g}_{2i} Q_2^{2i-1}) | n\rangle, \qquad m=1,2,... \cr}}

In the abovementioned case of the Ising model with only $g_2={\tilde g}_2=c$ and 
$g_4={\tilde g}_4=g$ non-zero,
we must take $Q_1=\sigma+\sigma^{-1}{\hat r}^{(1)}+\sigma^{-3}{\hat r}^{(2)}$ and
$Q_2=\sigma^{-1}{\hat v}+\sigma^{-1}{\hat x}^{(1)} \sigma^2+\sigma^{-1}{\hat x}^{(2)} \sigma^4$, and
eq.\eqdert\ reduces to
\eqn\geoising{\eqalign{
V_n&=1+c R_n +g R_n(R_{n+1}+R_n+R_{n-1})+g (R_n^{(2)}+R_{n+1}^{(2)}+R_{n+2}^{(2)})\cr
X_n^{(1)}&= c+ g(R_n+R_{n-1}+R_{n-2})\cr
X_n^{(2)}&= g \cr
R_n&=c V_n + g V_n( V_{n+1}X_{n+2}^{(1)}+V_n X_{n+1}^{(1)}+V_{n-1}X_n^{(1)}) \cr
R_n^{(2)}&= g V_n V_{n-1} V_{n-2}\cr}}
We first remark that $R_n=V_n X_{n+1}^{(1)}$ by comparing the second and fourth lines of
eq.\geoising. This is a particular case of a general duality
between $Q_1$ and $Q_2$ in the symmetric case when $g_i={\tilde g}_i$ for all $i$,
where we may write $Q_1=Q_2^\dagger$, where $\sigma^\dagger=\sigma^{-1} v$,
$(AB)^\dagger=B^\dagger A^{\dagger}$ for all operators $A$, $B$,
and $f^\dagger=f$ for all diagonal operators. Here this implies $\sigma^{-1}x^{(1)} \sigma^2=
r^{(1)} v^{-1} \sigma$ and $\sigma^{-1}x^{(2)} \sigma^4=r^{(2)}(v^{-1} \sigma)^3$,
i.e. $R_n=V_n X_{n+1}^{(1)}$ and $R_n^{(2)}=X_{n+1}^{(2)} V_n V_{n-1}V_{n-2}$.
Finally
eliminating $X_n^{(1)}$ and $R_n^{(2)}$ from eq.\geoising, we are left with 
\eqn\syising{\eqalign{
V_n(1-g^2(V_{n+1}V_{n+2}&+V_{n+1}V_{n-1}+V_{n-1}V_{n-2}))=1+R_n(c+g(R_{n+1}+R_n+R_{n-1})\cr
R_n&=V_n(c+g(R_{n+1}+R_n+R_{n-1})\cr}}

An important remark is in order about the generating function $R_n$. 
Although $n$ has the meaning of a maximal geodesic distance
between the two legs in the bipartite graph picture, it loses somewhat of its meaning
in the correspondence with Ising model configurations. Indeed, within configurations
of the Ising model on tetravalent planar graphs, $n$ is not the obvious geodesic distance
between the faces adjacent to the in- and out-coming legs, as its definition involves first
transforming the graph into a bipartite one, and it then corresponds to a distance where
edge-crossing is permitted only if the white vertex is on the right. 
This restriction is almost irrelevant, as there are in general sufficiently many
successions of black, white, black... bivalent vertices to allow for crossing
edges in both directions. One case however is troublesome: when an Ising edge 
connects a black vertex to a white one, in the absence of intermediate
bivalent vertices (not necessary here as the bicoloration is already ensured), 
the edge may only be crossed in one direction, leaving the white vertex on the right.
This introduces a bias in the notion of distance, having to do with the matter configurations
on the graph. 
An analogous situation was encountered in \NOUSHARD\ \GEOD\ in the case of hard dimers on
tetravalent planar graphs, namely on configurations of tetravalent planar graphs
where edges may (or may not) be occupied by dimers
which repel one-another in such a way that no two adjacent edges can be simultaneously 
occupied. In this case indeed, the notion of geodesic distance is biased by the dimers,
in that the occupied edges cannot be crossed in paths from $F_0$ to $F_1$. In both Ising and
hard-dimer cases,
the matter interfers with the space, by modifying the rules governing distances.
Consequently, we think it is interesting to investigate the dependence of the Ising two-leg
diagrams on this special distance, and it may eventually be that its difference with 
the true geodesic distance becomes irrelevant in large graphs.

\newsec{Exact solutions}

\subsec{Finding exact solutions: a general scheme}

All the equations listed in Sect. 5, despite their diversity, are all basically of the same form:
(possibly coupled) algebraic recursion relations expressing $R_n$ (and the other generating
functions involved) in terms of a finite number of previous terms $R_{n-1},R_{n-2},...,R_{n-k}$. In
principle the boundary data needed to entirely determine $R_n$ should consist of $k$ consecutive
initial values of $R_j$. It turns out however that we may drastically simplify, namely divide
by $2$ this required number of initial data by simply requiring that $\lim_{n\to \infty} R_n$
exists, and that it moreover coincides with the generating function $R$. Indeed, this is nothing
but restating the definitions of $R_n$ and $R$, as the latter was first obtained regardless
of the geodesic distance between legs, while the limit $n\to \infty$ of the former amounts to removing
the geodesic distance constraint in the counting of graphs. 
In the same fashion, all other generating functions involved tend to their obvious
limiting values when $n\to \infty$. This allows to linearize the various recursion
relations at large $n$, by setting say $R_n=R-\rho_n$, and similarly
for the other generating functions involved. At first order in $\rho_n$ and its other
counterparts, we obtain (possibly coupled) linear recursion relations. We immediately deduce
that $\rho_n \sim x^n$ (or a linear combination involving $x$'s of the same modulus) 
for some solution $x$ (with modulus less than 1) 
to the characteristic equation of the linearized recursion relations. 

To completely solve our equations, we start by determining the exact form of the
linearized solution, in general a linear combination $\rho_n\sim \sum_{j=1}^k a_j (x_j)^n$ where
$x_j$, $j=1,2,...,k$ denote all the (generically distinct)
solutions of the linearized characteristic equation with modulus less than 1.
In a second step, we obtain order by order in the $(x_j)^n$ the higher order contributions
to the true solution, expanded at large $n$. These take in general the form of recursion
relations for the coefficients of the multiple expansion in powers of the $(x_j)^n$. Solving these
recursion relations, and resumming the resulting series allows us to finally obtain
compact expressions for the exact solutions to the non-linear recursion relations at hand.
The result still depends on the initial parameters $a_j$, $j=1,2,...,k$, which then are fixed by
requiring that the terms involving the $k$ first $R_{-k}, R_{-k+1},...,R_{-1}$ drop out of 
the recursion relations.

In the following sections, we simply present the solutions, as we have found them. 
A case by case proof by substitution is left as an exercise to the reader. In many 
situations, the proof boils down to a certain identity between Chebyshev polynomials of the
first kind, as will be apparent soon. 

\subsec{Tetravalent case}

For pedagogical purposes, we detail in this simple case the general scheme presented in
the previous section.
Substituting $R_n=R-\rho_n$ into eq.\quaca\ we get at first order in $\rho_n$:
\eqn\quafirst{ \rho_n(1-3 gR)=gR(\rho_{n-1}+\rho_n+
\rho_{n+1})+O(\rho_n^2)} 
The linearized characteristic equation therefore reads
\eqn\licarqua{ 1-gR\left(x+{1\over x}+4\right) =0 }
For $g<g_c=1/12$, there is generically a unique solution $x\equiv x(g)$ with modulus
less than 1 to this equation, and we find that at first order $\rho_n=a x^n +O(x^{2n})$. 
We may now infer the general form $\rho_n=\sum_{j\geq 1} a_j x^{nj}$, $a_1=a$ for the
complete solution, where the coefficients $a_j$ are to be determined order by order in $x^n$.
We find explicitly
\eqn\wefiqua{ a_{k+1}=
\sum_{j=1}^k \left({x^j+{1\over x^j}+1\over x^{k+1}+{1\over x^{k+1}}-x-{1\over x} }
\right) a_j a_{k+1-j}}
solved recursively as
\eqn\solarho{a_k=a{1-x^k\over 1-x}\left({a x\over (1-x)(1-x^2)}\right)^{k-1}}
Picking $a=x(1-x)(1-x^2)\lambda$, $R_n$ is easily resummed
into
\eqn\fondaqua{ R_n = R {u_n u_{n+3}\over u_{n+1} u_{n+2}}, \qquad u_n=1-\lambda x^{n+1} }
This is the general solution to eq.\quaca, that converges for large $n$. 
To see why, it is simplest to substitute the form \fondaqua\ into the initial equation
\quaca\ which then boils down to the following quartic relation
\eqn\relaboilqua{ 
u_nu_{n+1}u_{n+2}u_{n+3}= {1\over R} u_{n+1}^2u_{n+2}^2+g R(u_{n-1}u_{n+2}^2u_{n+3}+
u_n^2u_{n+3}^2+u_n u_{n+1}^2u_{n+4})}
Substituting $u_n=1-\lambda x^{n+1}$ into this, we just have to check that the zeros
of the l.h.s. as a degree 4 polynomial in $\Lambda=\lambda x^n$ match those of the r.h.s.
as moreover the equation reduces to eq. \ralaRqua\ for $\Lambda=0$.

Finally, the ``integration"
constant $\lambda$ is now fixed by further requiring that eq.\quaca\ makes sense at $n=0$,
in which case the term $R_{-1}$ must drop off the r.h.s. of the recursion relation,
in other words we have to impose $R_{-1}=0$. This simply gives $\lambda=1$, and finally the
exact solution to our combinatorial problem reads
\eqn\exacombqua{ R_n=R {(1-x^{n+1})(1-x^{n+4})\over (1-x^{n+2})(1-x^{n+3})}=
R {U_n U_{n+3}\over U_{n+1}U_{n+2} }, \qquad 
U_n\equiv U_n\left(\sqrt{x}+{1\over \sqrt{x}}\right)={x^{n+1\over 2}-x^{-{n+1\over 2}}\over
x^{1\over 2}-x^{-{1\over 2}}}}
where we have identified the Chebyshev polynomials $U_n$ of the first kind.

As a simple application, the formula \exacombqua\ gives access to the generating function
for two-leg diagrams whose legs lie in the same face, already identified as that
of rooted tetravalent planar maps or quadrangulations. We find that
\eqn\rzerqua{ R_0= R {x+{1\over x} \over x+{1\over x}+1}=R {1-4gR\over 1-3gR}
=R-gR^3}
with $R$ as in \Rfqua. This
result was first obtained by Tutte \TUT\ in a completely different, though combinatorial, manner.

\subsec{Even valences}

Let us for definiteness consider the equation \arbeven\ with only $g_4,g_6,...,g_{2m+2}$
non-zero. Linearizing again the equation at large $n$ and solving for the leading $\rho_n\sim
x^n$, we find that $x$ must obey the following characteristic equation:
\eqn\chargen{\eqalign{0=\chi_m(x)&\equiv 1 -\sum_{k=0}^m  g_{2k+2} R^{k}
\sum_{l=0}^{k} {2k+1\choose l} {1\over x^{k-l}} {1-x^{2k+1-2l}\over 1-x}\cr
&=1 -\sum_{k=1}^m  g_{2k+2} R^{k} \sum_{l=0}^{k} {2k+1\choose l} U_{2k-2l}(w) \cr}}
expressed in terms of Chebyshev polynomials of $w=\sqrt{x}+1/\sqrt{x}$.
Note that $\chi_m(x)$ is actually a degree $m$ polynomial in $x+1/x$, with
generically $m$ distinct roots with modulus less than 1, denoted by $x_1,x_2,...,x_m$.
Repeating the straightforward, though tedious, exercise of previous section,
we end up with the following exact solution
\eqn\mxsol{\eqalign{
&R_n=R{u_n^{(m)} u_{n+3}^{(m)} \over u_{n+1}^{(m)} u_{n+2}^{(m)} } \cr
&u_n^{(m)}= \sum_{l=0}^m (-1)^l \sum_{1\leq m_1 <...<m_l\leq m}
\prod_{i=1}^l \lambda_{m_i} x_{m_i}^{n+m} \prod_{1\leq i<j\leq l} c_{m_i,m_j}\cr
&c_{a,b}\equiv {(x_{a} -x_{b})^2 \over (1-x_{a} x_{b})^2} \cr
}}
Remarkably, the structure of $u_n^{(m)}$ matches exactly that of an $N$-soliton
tau-function of the KP hierarchy \JM\ which reads
\eqn\statjimi{\eqalign{
&\tau=\sum_{r=0}^N \sum_{i_1<...<i_r} \prod_{\mu=1}^r e^{\eta_{i_\mu}}
\left(\prod_{\mu<\nu} c_{i_\mu,i_\nu}\right) \cr
&c_{i,j}={(p_i-p_j)(q_i-q_j)\over (p_i-q_j)(q_i-p_j)}\cr}}
where the $\eta$'s contain the times' dependence of the KP hierarchy. 
Our solution \mxsol\ simply amounts to identifying $N=m$, $e^{\eta_i}=-x_i^{n+m} \lambda_i$,
$p_i=x_i$ and $q_i=1/x_i$, $i=1,2,...,m$. 
This surprising relation suggests the existence of an underlying integrable structure
for our recursion relations, but remains mysterious. 

Again, to further fix the ``integration constants" $\lambda_1,\lambda_2,...,\lambda_m$,
we simply have to express that $u_{-1}=u_{-2}=...=u_{-m}=0$, with the result
\eqn\resvani{\lambda_i =\prod_{j\neq i} {1-x_i x_j \over x_i-x_j} \qquad i=1,2,...,m }
and finally the complete solution to \arbeven\ reads
\eqn\solarb{ R_n=R {U_n(w_1,...,w_m) U_{n+3}(w_1,...,w_m) \over
U_{n+1}(w_1,...,w_m) U_{n+2}(w_1,...,w_m)} }
where 
\eqn\detru{ U_n(w_1,...,w_m)\equiv
\det\left[ U_{n+2j-2}(w_i) \right]_{1\leq i,j \leq m}}
in terms of Chebyshev polynomials of the first kind expressed at the values
$w_i=\sqrt{x_i}+1/\sqrt{x_i}$, $i=1,2,...,m$.
The precise proof of these statements can be found in \GEOD.

For illustration, in the tetra/hexavalent case $m=2$, we have the characteristic
equation
\eqn\tetrhexchar{
0 =\chi_2(x)\equiv 1-g_4 R\left(x+{1\over x}+4\right)-g_6 R^2
\left(x^2+6 x+16+{6\over x}+{1\over x^2}\right) }
and for instance the result \solarb\-\detru\ reads for $n=0$:
\eqn\initdim{ R_0= R{1+(x_1+{1\over x_1})(x_2+{1\over x_2})\over
1+(1+x_1+{1\over x_1})(1+x_2+{1\over x_2})}=R{1-4g_4R-15g_6R^2\over
1-3 g_4R-10 g_6 R^2}}
where we have reexpressed the symmetric functions of $x_1+1/x_1$ and $x_2+1/x_2$
in terms of the coefficients of $\chi_2(x)$.
The result is nothing but the generating function for rooted tetra/hexavalent 
planar graphs.

\subsec{Trivalent case}

Repeating the exercise of Sect.6.2 for the trivalent case of eqs. \transfortri, 
we have found the following solution
\eqn\foundtri{ R_n=R {u_n u_{n+2} \over u_{n+1}^2}, \qquad u_n=1-\lambda x^{n+1}}
while
\eqn\ssoltri{ S_n=S - {t_n \over u_{n}u_{n+1}}, \qquad 
t_n=gR^2(1-x)(1-x^2)\lambda x^n} 
In eqs.\foundtri\-\ssoltri, the parameter $x$ is the unique solution to
the linearized characteristic equation
\eqn\linchartri{1-g^2R^3\left(x+{1\over x}+2\right)=0}
such that $|x|<1$.

Requiring that $R_{-1}=0$, fixes $\lambda=1$, and we obtain  
the complete solution
\eqn\completri{\eqalign{
R_n&=R {(1-x^{n+1})(1-x^{n+3})\over (1-x^{n+2})^2} \cr
S_n&=S -gR^2{ (1-x)(1-x^2) x^n \over (1-x^{n+1})(1-x^{n+2})}\cr}}
from which we read off the following compact expressions
for $R_0$ and $S_0$, respectively generating two- and one-leg diagrams
of trivalent planar graphs, with respectively the two legs in the same face
and the leg in the external face:  
\eqn\valzertri{\eqalign{
R_0&= R {x+{1\over x}+1\over x+{1\over x}+2} =R-g^2 R^4 \cr
S_0&=S-gR^2\cr}}
In particular, $R_0$ is the generating function for rooted trivalent planar graphs.

\subsec{Arbitrary valences}

The case of arbitrary valences still awaits a good solution, in the same spirit as 
the case of even valences. It is however possible to find ``integrable-like" solutions
to the corresponding recursion relations (containing only one integration
constant), which for the time being are still too restrictive
to describe the general case.
Indeed, we need a sufficient number of integration constants
to allow for satisfying all the necessary initial conditions of our combinatorial problem.
This number is exactly the degree of the characteristic equation of the linearized recursions,
when expressed as a polynomial in $x+1/x$.

For arbitrary values of $g_3,g_4,g_5,...$ we have the following 
``one-$x$" solutions:
\eqn\soltriter{\eqalign{
R_n&= R {(1-\lambda x^{n+1})(1-\lambda x^{n+3})\over (1-\lambda x^{n+2})^2} \cr
S_n&=S -\sqrt{R x}\, {(1-x)^2\lambda x^n\over
(1-\lambda x^{n+1})(1-\lambda x^{n+2})} \cr}}
where $x$ is any solution with modulus less than 1 to the corresponding linearized characteristic
equation. The latter reads for instance in the case of tri/tetravalent graphs
where only $g_3$ and $g_4$ are non-zero:
\eqn\lintritet{\left(g_4 R(x+{1\over x}+4)+S(2g_3+3g_4S)-1\right)^2-R(g_3+3 g_4 S)^2(x+{1\over
x}+2)=0}
and amounts to 
\eqn\amto{ \sqrt{R} (g_3+3 g_4 S)(\sqrt{x}+{1\over \sqrt{x}})=1-g_4 R(x+{1\over x}+4)
-S(2g_3+3g_4S) }
as $R$ has a power series expansion of the form $R=1+O(g_3,g_4)$.
In this particular case, it would be desirable to obtain the full solution
involving the two roots $x_1$ and $x_2$ of \lintritet\ and two integration constants, 
to be able to solve simultaneously the two boundary conditions $R_{-1}=0$ and $R_{-1}S_{-1}=0$
(to be understood as $\lim_{n\to -1} R_n S_n=0$)
obtained from \tritetr\ at $n=0$, and clearly not satisfied by \soltriter.

\subsec{Bipartite $p$-valent case}

Repeating the usual exercise with eq.\onlybw\ in which we first
eliminate $X_n$, we have found the following solution, including
only one integration constant:
\eqn\solpurepcolonex{ R_n=R { u_nu_{n+p+1}\over u_{n+1} u_{n+p}}, \qquad u_n=1-\lambda x^{n+1} }
valid provided $x$ is chosen among the roots of modulus less than 1 of the linearized
characteristic equation:
\eqn\lincarpcolonex{1- g {\tilde g}_1 R^{p-2} {1\over x^{p-2}}(1+x+x^2+...+x^{p-2})^2 =0}
As noticed in the case of planar graphs of arbitrary valences, the degree of this
polynomial of $x+1/x$ is however $p-2$, hence only for $p=3$ is the solution \solpurepcolonex\
completely general. This is the case of bipartite trivalent planar graphs, for which eq.\onlybw\
reduces to $R_n=1+g{\tilde g}_1 R_n(R_{n+1}+R_{n-1})$.
In this case, the solution is further fixed by requiring $R_{-1}=0$, hence $\lambda=1$,
and it reads
\eqn\tricase{ R_n=R {(1-x^{n+1})(1-x^{n+5})\over (1-x^{n+2})(1-x^{n+4})}}
with $1-g{\tilde g}_1 R (x+1/x+2)=0$, $|x|<1$.
$R_n$ is the generating function for bipartite trivalent graphs with two legs attached
to vertices of opposite colors, and geodesic distance between those less or equal to $n$. 
In particular, 
\eqn\rzereultri{ R_0=R {x^2+{1\over x^2}+x+{1\over x}+1\over (x+{1\over x})(x+{1\over
x}+2)}=R{1-3g {\tilde g}_1R+g^2{\tilde g}_1^2R^2\over 1-2g {\tilde g}_1 R}}
is the generating function for rooted bipartite trivalent planar graphs (with weights
$g$/${\tilde g}_1$ per black/white vertex)

For $p\geq 4$, we must work out the generalizations of \mxsol, which now read as follows.
Introducing
\eqn\intropqc{ \eqalign{ p(x)&=x(1+x+...+x^{p-2}), \qquad q(x)= p(1/x) \cr
p_i&=p(x_i), \qquad q_i=q(x_i) \cr
c_{i,j}&={(p_i-p_j)(q_i-q_j)\over (p_i-q_j)(q_i-p_j)} \cr}}
where $x_i$, $i=1,2,...,p-2$ denote the generically distinct roots of eq.\lincarpcolonex\
with modulus less than 1, the general solution now takes the form
\eqn\forpcolmax{R_n=R {u_n^{(m)} u_{n+p+1}^{(m)}\over u_{n+1}^{(m)}u_{n+p}^{(m)}} }
with $m=p-2$ and as before
\eqn\kxun{
u_n^{(m)}=\sum_{l=0}^m (-1)^l \sum_{1\leq i_1<i_2<...<i_l\leq m} 
\prod_{t=1}^l \lambda_{i_l} x_{i_l}^n
\prod_{1\leq r<s\leq l} c_{i_r,i_s}}
Note the absolutely remarkable fact that we obtain again an expression in terms of
the tau-function for the KP hierarchy, but with different kinematics, in the form
of an implicit relation between the $p$'s and $q$'s of eq.\intropqc.

Imposing moreover the vanishing of the first terms $u_{-1}=u_{-2}=...=u_{-p+1}$,
we find that
\eqn\valulam{
\lambda_i= x_i^{(p-1)m-1} \prod_{j\neq i} {q_i-p_j\over p_i-p_j} }
This fixes completely the solution to our combinatorial problem, and in particular
gives a compact expression for the generating function $R_0$ for rooted
bipartite $p$-valent planar graphs. For $p=4$ for instance, it reads
\eqn\pfourcol{
R_0=R{1-5g {\tilde g}_1R^2+3g^2 {\tilde g}_1^2 R^4\over 1-3 g {\tilde g}_1R^2}}
while $R$ satisfies $R=1+3 g {\tilde g}_1 R^3$.

\subsec{Constellations}

The solution for general $p$-constellations is similar to that for pure $p$-valent bipartite
graphs. Indeed, we find that the general solution has the exact same form
\forpcolmax\ as in the previous section, with $u_n^{(m)}$ given by \kxun\ and 
$p_i,q_i,c_{i,j}$ defined as in \intropqc. The only difference is that now $m$ may take
a larger value say $(p-1)k-1$ in the case when only $g,{\tilde g}_1,{\tilde g}_2,...,{\tilde
g}_k$ are non-zero (namely of constellations with white $p$ valent vertices, and black
vertices with valences $p$, $2p$, ...,$kp$). 
The $x$'s entering these formulas are now the generically distinct roots
with modulus less than 1 of the linearized characteristic equation, itself a polynomial
of degree $(p-1)k-1$ of $x+1/x$, reading
\eqn\carpconstel{
1=\left(1+{1\over x}+...+{1\over x^{p-2}}\right)
\sum_{i\geq 1} g^i {\tilde g}_i R^{(p-1)i-2}
\sum_{m=0}^{(p-1)i-1} \sum_{j=0}^{i-1} {j+m\choose m}{pi-2-j-m\choose i-j-1}x^{m-j(p-1)}}

For illustration, in the case of $3$-constellations with say only $g$ and 
${\tilde g}_2$
non-zero ($R_n$ obeys eq.\instcons\ with ${\tilde g}_1=0$), 
we get the linearized characteristic equation:
\eqn\charthree{1-g^2 {\tilde g}_2R^3(x+{1\over x}+2)(x^2+{1\over x^2}+x+{1\over x}+6)=0 }
of degree $3$ in $x+1/x$. Picking the three roots with $|x_i|<1$,
we finally get the generating function for rooted 3-constellations with white trivalent
and black hexavalent vertices
\eqn\extricons{ R_0=R{1-17 g^2 {\tilde g}_2R^3+25 g^4 {\tilde g}_2^2 R^6\over 1-10
g^2 {\tilde g}_2 R^3} }
while $R$ satisfies the equation \singal\ with only $g$ and ${\tilde g}_2$ non-zero, namely
$R=1+10 g {\tilde g}_2 R^4$.

\subsec{Ising model}

We have not been able to find a nice structure for the general solution of the even-valent
bipartite case in general. 
In the particular case of the Ising model
with zero magnetic field (eq.\syising), we have been able to derive all possible solutions
involving only one integration constant. 

The usual linearization of eq.\syising\ yields the following characteristic equation
\eqn\carisi{\eqalign{
\left(x+{1\over x}+{c\over gV (1-3gV)}-{1-gV\over gV}\right)&
\left(x+{1\over x}-{c\over gV (1-3gV)}-{1-gV\over gV}\right) \times \cr
&\times \left(x+{1\over x}+{1+gV\over gV}\right) =0\cr}}
of degree $3$ in $x+1/x$. As already discussed before, we would need in principle 
to find the full solution
to \syising\ that converges at large $n$, including all three $x_1,x_2,x_3$, respectively
the roots of the three factors in
eq.\carisi\ with modulus less than 1, and therefore including also 
three integration constants. We now display the 
solutions with one $x$ for each of the three factors in \carisi.

For both $x=x_1$ and $x=x_2$, we have found that
\eqn\solisione{ V_n=V{u_n u_{n+3}\over u_{n+1} u_{n+2}},\qquad u_n=1-\lambda x^n}
while for $x=x_1$:
\eqn\isir{ R_n=R- V {\lambda (1-x)(1-x^2) x^n\over u_{n+1}u_{n+2} } }
and for $x=x_2$:
\eqn\isirtwo{ R_n=
R+V {\lambda (1-x)(1-x^2) x^n\over u_{n+1}u_{n+2} } }
For $x=x_3$ however, the solution is quite different:
\eqn\solisithree{\eqalign{
V_n&=V {u_n u_{n+3}\over u_{n+1} u_{n+2}}\cr 
u_n&=1-2 \lambda x^n-z \lambda^2 x^{2n} \cr
z&={\big((c-2)(x+{1\over x})+c-8\big)\big((c+2)(x+{1\over x})+c+8\big)
\over \big(x^2+{1\over x^2}-(c-4)(x+{1\over x})-(c-2)\big)
\big(x^2+{1\over x^2}+(c+4)(x+{1\over x})+(c+2)\big)}\cr}}
The growing complexity of the solutions 
lets us expect a quite involved general three $x$-solution, yet to be found.

\newsec{Continuum limit}

The exact solutions of Sect. 6 allow us to investigate the properties of the corresponding
classes of planar graphs in terms of the geodesic distance, in the limit
when the latter becomes large. 
This limit must clearly be taken simultaneously with the so-called critical
limit of large graphs, reached in turn
by letting the various weights per vertex approach some critical locus, corresponding to
approaching the finite radii of convergence of the various combinatorial series involved.

The large $n$ limit of $R_n=R+ a_1 x_1^n+...$ where $x_1$ is the largest root with modulus
less than 1 of
the characteristic equation, leads to the natural definition
of a correlation length $\xi=-{\rm Log}\, |x_1|$, governing the exponential decay
of $R_n-R\sim \exp(-n/\xi)$ as a function of the geodesic distance $n$.
A good continuum limit may therefore be reached by letting $x_1$ (and possibly other $x$'s)
tend
to the value $1$, while simultaneously keeping $n/\xi$, the continuum geodesic
distance, fixed. This in turn implies certain
relations between the vertex weights are reached, via the characteristic equation relating
them to the $x$'s. These relations express nothing but the abovementioned critical limit, which must
therefore be taken simultaneously with the continuum one. 

\subsec{Tetravalent case}

We illustrate the above with the case of tetravalent graphs, with $R_n$ given by \exacombqua.
The critical limit $g\to g_c=1/12$ is reached by taking say
\eqn\critlim{ g=g_c(1-\epsilon^4), \qquad R={R_c \over 1+\epsilon^2} }
with $R_c=2$, and the solution of eq.\licarqua\ with modulus less than 1 reads 
\eqn\xlim{ x(g)= {1+2\epsilon^2-\epsilon\sqrt{3(2+\epsilon^2)} \over 1-\epsilon^2} }
hence $x=e^{\epsilon\sqrt{6}}+O(\epsilon^2)$ as $\epsilon\to 0$, and $\xi\propto 1/\epsilon$. 
Finally setting 
\eqn\scalged{n={r\over \epsilon}}
we may simply express the continuum limit
of the quantity $R_n$, or more interestingly that of $R-R_n$,
generating two-leg diagrams of tetravalent planar graphs with distance at least
$n$ between the two legs.
We find that
\eqn\contqua{{\cal F}(r)=\lim_{\epsilon\to 0}{R-R_n\over \epsilon^2 R}
={3\over \sinh^2\left(\sqrt{3\over 2} r\right)} }
This is the continuum two-point correlation function for random surfaces with two marked
points at (rescaled) geodesic distance larger or equal to $r$.
It coincides with the scaling function derived in \AW, by use of transfer matrix formalism.
This in turn yields  
the continuum two-point correlation function for random surfaces with two marked
points at (rescaled) geodesic distance $r$:
\eqn\eqr{{\cal G}(r)=-{\cal F}'(r)= 3\sqrt{6} {\cosh\left(\sqrt{3\over 2}r\right)\over
\sinh^3\left(\sqrt{3\over 2}r\right) } }

This result may in turn be interpreted in terms of graphs with large but finite number
$N$ of vertices.
Indeed, the above scaling relations \critlim\ and \scalged\ 
imply the following relation between the correlation length and the deviation from
the critical point $\xi \sim (g_c-g)^{-\nu}$, with the exponent $\nu=1/4$, and
therefore the fractal dimension $d_F=1/\nu=4$ for the present model of random
planar surfaces. More concretely, this tells us in particular that the number of
faces in a large graph lying at geodesic distance $\leq n$ from the external one 
behaves as
$n^{d_F}=n^4$. A simple measure of this number is indeed the ratio $R_n\vert_{g^N}/R_0\vert_{g^N}$
of the corresponding coefficients of $g^N$ in the two power series expansions, giving
the proportion of graphs with the two legs distant by at most $n$ to that with the two legs
in the same face. Using our exact solution and performing a saddle-point expansion,
we find
\eqn\measureofdf{ \lim_{N\to \infty} {R_n\vert_{g^N} \over R_0\vert_{g^N}}\sim {3\over 56} n^4}
giving an explicit illustration of the fractal dimension 4.

This suggests to set $n=r N^{1/4}$ and to write 
$R_n\vert_{g^N}$, again by use of a saddle-point expansion, as
\eqn\asymptorn{R_n\vert_{g^N}\sim {4\over \pi} {(12)^N\over N^{3/2}} \int_0^\infty
du u^2 \left(1+{\rm Re}\, {\cal F}(r\sqrt{-iu})\right) }
and similarly for $R\vert_{g^N}$ with $\cal F$ replaced by $0$. The ratio
$R_n\vert_{g^N}/R\vert_{g^N}$ gives the probability $P(r)$  
for a random surface with two marked points that their geodesic distance be less or equal to $r$:
\eqn\probar{P(r)={2\over \sqrt{\pi}}\int_0^\infty du u^2 e^{-u^2} \left(1+{\rm Re}\, {\cal
F}(r\sqrt{-iu})\right) }

\subsec{Arbitrary even valences, multicriticality}

Repeating the analysis of previous section for the case of graphs with even valences say up
to $2m+2$, various critical points may be reached in the space of weights $g_j$. 
Concretely, picking the particular weights that ensure the following form for \relgeR:
\eqn\newrelgeR{ {g_c-g\over g_c}= \left({V_c-V\over V_c}\right)^{m+1} }
where $V=g R$, $g_4=g$, $g_{2j}= g^{j} z_{j}$ fixed by the form \newrelgeR\ and
$V_c=m/6$, $g_c=m/(6(m+1))$,
the corresponding multicritical limit is obtained by setting
\eqn\multricri{ g=g_c(1-\epsilon^{2(m+1)}),\qquad R=R_c {1-\epsilon^2\over 1-\epsilon^{2(m+1)}}}
with $R_c=V_c/g_c=m+1$. Remarkably, the characteristic equation \chargen\
then turns into
\eqn\turcar{ \chi_m(x)= \left({V_c-V\over V_c} \right)^m P_m\left(
{1-\epsilon^2\over \epsilon^2} \big(x+{1\over x}-2\big)\right) =0 }
for some fixed degree $m$ polynomial 
\eqn\pmpol{P_m(u)=\sum_{l=0}^m  (-u)^l {l! \over (2l+1)!} {m! \over (m-l)!}}
This allows to get
the leading behavior of the various $x$'s as
\eqn\varbeha{ x_i=e^{-a_i \epsilon} +O(\epsilon^2)}
where $a_i^2$ are the roots of $P_m$, and $a_i$ are taken with positive real part.
Consequently all the $x$'s tend to 1 simultaneously, and the correlation length of the problem
reads $\xi\sim 1/\epsilon$, hence we now have the relation $\xi\sim (g_c-g)^{-\nu}$
with $\nu=2(m+1)$, i.e. a fractal dimension $d_F=2(m+1)$.

The multicritical continuum limit is therefore still obtained by setting \scalged, and
the exact solution \mxsol\ yields the following two-point correlation of random
surfaces with multicritical weights (known to simulate non-unitary matter conformal field theories
with central charges $c(2,2m+1)=1-3(2m-1)^2/(2m+1)$ coupled to two-dimensional quantum gravity), 
with two marked points at geodesic distance less or equal to $r$:
\eqn\genscafunc{ {\cal F}(r)=-2 {d^2\over dr^2} {\rm Log}\,
{\cal W}\left(\sinh\left(a_1{r\over 2}\right), \sinh\left(a_2{r\over 2}\right),...,
\sinh\left(a_m{r\over 2}\right)\right) }
where ${\cal W}(f_1,f_2,...,f_m)$ stands for the Wronskian determinant $\det\big[ f_i^{(j-1)}
\big]_{1\leq i,j \leq m}$.

It is known from matrix model solutions that the general case of arbitrary valences leads 
to the same multicritical points. In particular we expect the scaling
functions \genscafunc\ to be the same at these points. The same remark applies to constellations
as well. New multicritical points corresponding to conformal theories with central charges
$c(p,q)=1-6(p-q)^2/(pq)$ for $p,q$ two coprime integers can be reached within the framework
of two-matrix models, corresponding to the general bipartite graphs. In the latter case, we expect
some new scaling functions, characteristic of these other universality classes.
An example will be given in next section, when discussing the Ising model.

\subsec{Critical/continuum limit in general}

The (multi-) critical continuum limits are reached by letting a number of the
roots $x$ of the characteristic equation tend to 1 simultaneously.  

An alternative route for deriving the critical continuum limit of say the tetravalent case
would have been to postulate the form $R_n=R(1-\epsilon^2{\cal F}(n\epsilon))$ and
plug this ansatz into the recursion relation \quaca. With $g=g_c(1-\epsilon^4)$, expanding
the equation up to order $4$ in $\epsilon$, we arrive at the following differential equation
for ${\cal F}$: 
\eqn\painleve{{\cal F}''(r)-3{\cal F}^2(r)-6{\cal F}(r)=0}
The function $\cal F$ of eq.\contqua\ is the unique solution to \painleve\
such that ${\cal F}(0)=\infty$ and ${\cal F}(\infty)=0$. 

More generally, we may derive differential equations for the multicritical cases of Sect.7.2
as well. These take the form
\eqn\fordif{ {\cal R}_{m+1}[1+{\cal F}] ={\cal R}_{m+1}[1] }
where ${\cal R}_m[u]$ is the $m$th KdV residue $(d^2-u)^{m-1/2}\vert_{d^{-1}}$
where $d\equiv d/dr$ \GD. For instance, in the case of multicritical tetra/hexavalent graphs
($m=2$) eq.\fordif\ reads
\eqn\painlesix{ {\cal F}^{(4)}(r) -10{\cal F}(r) {\cal F}''(r)-10 {\cal F}''(r)
-5 ({\cal F}'(r))^2+10
({\cal F}(r))^3 +30({\cal F}(r))^2 +30{\cal F}(r) =0}
 
Even more generally, recall that the KdV residues naturally arise
(e.g. in the context of matrix model solutions) when solving
the differential operator equation $[P,Q]=1$, where $Q=d^2-u$ and $P$ some 
degree $2m+1$ differential operator. The equation indeed boils down to
$2d/dr({\cal R}_{m+1}[u])=1$. In the present case, we rather have to write the equation $[P,Q]=0$,
which turns into ${\cal R}_{m+1}[u]=$const.

At the discrete level, comparing say \quaca\ with the equations determining the
matrix model solution for tetravalent graphs of arbitrary genus, we simply would
have to replace $1$ by $n/N$ in the r.h.s. of \quaca. 
We may think of the differential operator $Q=d^2-u$ as the continuum limit
of the operator $Q$ \opQ\ defined in Sect.5.2, now acting on functions
of the variable $r=n\epsilon$. 
In the other cases described in Sect.6, we always have such an operator $Q$ at hand  
and again the recursion relations resemble strongly those obtained in the solutions
of the multi-matrix models describing two-dimensional quantum gravity coupled with matter,
up to the same substitution $1\to n/N$. This suggests that the relevant (coupled)
differential equations governing the (multi-) critical continuum limit should read
$[P,Q]=0$, with $Q$ the continuum limit of the operator say $Q_1$ used in our approach,
taking the form of a differential operator of degree $q$ say, and $P$ a differential operator
of degree $p$ coprime with $q$. Our claim is that the generalized two-point functions
$\cal F(r)$ for random surfaces in the presence of critical matter (corresponding to
conformal field theories with central charges $c(p,q)<1$), with two marked points
at geodesic distance less than $r$, should be governed by $[P,Q]=0$, where 
$Q=d^q- qu d^{q-2}_...$, and $u=1+{\cal F}$. This is illustrated in the case of the Ising model in
the next section.

\subsec{Ising model}

The multicritical limit of the Ising model is obtained as follows.
Starting from the equations \finroneis\ rewritten
as $W(R)=0$, $W$ a polynomial of degree $5$, the
tri-critical points are determined by setting $W'(R)=W''(R)=0$, and we find
that $c_c=\pm 4$ while $g_c=10/9$, and $R_c=-3/5$, $V_c=-3/10$.
The tricritical limit is approached by setting 
\eqn\critli{ c=4, \quad R=R_c(1-\epsilon^2)\  \ \Rightarrow \ \ g=g_c(1-{16\over 5}\epsilon^6)}
We note that in the characteristic equation \carisi\ only the first and last factor
tend to $x+1/x-2$ as $\epsilon\to 0$, while the middle one tends to $x+1/x+10$.  
This means that only two of the three $x$'s tend to 1 in this limit.
More precisely, we have $x_3=e^{-\sqrt{6}\epsilon}+O(\epsilon^2)$ and 
$x_1=e^{-2\sqrt{3}\epsilon}+O(\epsilon^2)$.
This displays the fractal dimension of graphs with critical Ising configurations,
namely $d_F=6$, obtained by expressing the correlation length $\xi\sim (g_c-g)^{-\nu}$,
where $\nu=1/6$.
Recall however that the distance $n$ or its rescaled version $r$ are slightly different
from the true geodesic distance in the Ising tetravalent graphs, and the fractal dimension measured here
might be different from that associated to the true geodesic distance.

Further substituting 
\eqn\ansatisi{\eqalign{ R_n&=R(1-\epsilon^2{\cal F}(n\epsilon)), \qquad r=n\epsilon \cr
V_n&={R_n\over c+g(R_{n+1}+R_n+R_{n-1})} \cr}}
into the recursion relations \syising\ together with \critli, 
we find by expanding up to order 6 in $\epsilon$ that the two-point function ${\cal F}$
obeys the differential equation
\eqn\isingdif{ {\cal F}^{(4)}(r) -18{\cal F}(r) {\cal F}''(r)-18 {\cal F}''(r)
-9 ({\cal F}'(r))^2+24
({\cal F}(r))^3 +72({\cal F}(r))^2 +72{\cal F}(r) =0}
This equation is precisely what one would get by writing $[P,Q]=0$ for differential
operators $Q=d^3-3u d -3u'/2$ and $P$ of order $4$, $u=1+{\cal F}$.

Looking for convergent solutions in the form ${\cal F}(r)=a e^{-k r}$, we
find that $k_1=\sqrt{6}$ or $k_2=2\sqrt{3}$, corresponding to $e^{-k_1 r}=x_3^n$ and 
$e^{-k_2r}=x_1^n$ respectively. 
Proceeding like in the discrete case, we may now solve the differential
equation order by order in $e^{-k_ir}$, with a double power series expansion
${\cal F}(r)=\sum_{m,p\geq 0} a_{m,p} e^{-r(mk_1+pk_2)}$, in terms
of the two integration constants $\lambda=a_{1,0}$, $\mu=a_{0,1}$, while $a_{0,0}=0$. 
The differential
equation \isingdif\ indeed simply amounts to a recursion relation on the coefficients $a_{m,p}$.
Resumming the series for ${\cal F}(r)$ finally yields:
\eqn\iscafuncising{\eqalign{
{\cal F}(r)&= -{d^2\over dr^2}{\rm Log}\left(
1-{\lambda\over 6}e^{-r\sqrt{6}}
-{\mu\over 12}e^{-2r\sqrt{3}}-{\lambda^2\over 288}e^{-2r\sqrt{6}} \right.\cr
&\left.-{17-12\sqrt{2}\over 72}\lambda \mu e^{-r(\sqrt{6}+2\sqrt{3})} 
+{577-408\sqrt{2}\over 3456}\lambda^2\mu e^{-2r(\sqrt{6}+\sqrt{3})} \right)\cr} }
The two integration constants are further fixed by the boundary conditions.
The latter are obtained by requiring that the recursion relations \syising\
also make sense at $n=0,1$, namely $R_{-1}=V_{-1}=0$, while $\lim_{n\to 0} V_{n-1}V_{n-2}=0$. 
The first conditions give ${\cal F}(0)=\infty$, while the latter implies $R_{-1}R_{-2}=0$.
As in the case of planar graphs with only tetra/hexavalent vertices, this implies a
higher order vanishing of the argument of the logarithm in \iscafuncising, which plays the role
of continuum limit of $u_n$, while the condition implies that both $u_{-1}$ and $u_{-2}$ vanish.
This is easily solved into $\lambda=-12(17+12\sqrt{2})$ and $\mu=12(4+3\sqrt{2})$, so that finally
\eqn\calfisi{
{\cal F}(r)=-{d^2\over dr^2}{\rm Log}\Big(
\sinh\big(r(\sqrt{6}+\sqrt{3})\big)
+(17+12\sqrt{2})\sinh\big(r(\sqrt{6}-\sqrt{3})\big)
-2(4+3\sqrt{2})\sinh(r\sqrt{3})\Big)} 
and we also get the correlation ${\cal G}(r)=-{\cal F}'(r)$ for the Ising model
on random surfaces, with two marked points at (special) geodesic distance $r$. 

\newsec{Conclusion}

In this note, we have addressed the problem of enumeration of various types of
planar graphs with external legs, while keeping track of a suitable geodesic distance
between these legs. 
The basic tool we used are blossom-trees, namely trees carrying the minimal information
needed to close them back into planar graphs. This information is essentially
contained in the two types of leaves (black and white), and we have devised a
compact algebraic way of keeping track of the geodesic distance between legs,
by introducing operators $Q$ describing the structure of the rooted trees
around their first vertex.
This operator acts formally on a basis $|n\rangle$ indexed by relative integers,
and may as well be viewed as acting on sequences $\{p_n\}_{n\in \IZ}$, via the shift
operator $\sigma$ and its relative integer powers, and a number of (combinatorial)
diagonal operators.
This definition is clearly borrowed from that of the $Q$ operator of matrix models,
that generate the multiplication by an eigenvalue $\lambda$ on the basis of 
(right) monic (bi-)orthogonal polynomials $p_n(\lambda)$, $n=0,1,2...$
In this framework, the main recursion relations for the coefficients of $Q$
are obtained by considering the operator $P$, acting on the $p_n(\lambda)$
by  differentiation w.r.t. $\lambda$. The canonical relation $[P,Q]=1$
determines in fact all the functions of the problem, and turns into differential equations
in the critical scaling limit.
Here, by analogy with this case, we have been led to set $[P,Q]=0$, for some operator $P$
still awaiting a good combinatorial meaning. Nevertheless, this equation also turns into differential
equations for the physical quantities of our problem.
The mysterious part of this correspondence is that the natural  
variable in the matrix model approach is the rescaled cosmological constant, which generates
the topological expansion of the free energy, namely the expansion in powers of $r$ of the
function $u(r)$ has coefficients corresponding to planar graphs of fixed genus.
This must be contrasted with the present findings, where $r$ has the meaning of geodesic
distance.  This seems to indicate that a more general structure should exist, that includes
both the topological and geodesic directions, probably some suitably defined matrix model
of sorts.

In view of these strong analogies with matrix model solutions, we may want to characterize the class
of possibly decorated planar graphs for which a tree formulation exists as that for which
there exists a matrix model formulation admitting a solution via orthogonal polynomials.
This would exclude for instance the case of the three-state Potts model, a generalization of the 
Ising model whose configurations are graphs with vertices of three possible colors, and
edge weights $w_{a,b}=e^{K\delta_{a,b}}$ according to the colors $a,b$ of the adjacent vertices.
More generally, loop models on graphs also correspond to matrix models without orthogonal polynomial
solutions: their configurations are simply mutually- and self-avoiding loops drawn on the edges of 
planar graphs, with a weight $n$ per loop (the so-called O(n) model, extensively solved
in \EK). It would be extremely interesting to investigate any of these models using tree techniques.

A final striking outcome of our work is the emergence of soliton-like tau functions entering the 
explicit exact formulas  for a number of generating functions for two-leg diagrams
with legs distant by at most $n$. 
A reason why this should happen in the first place may be perhaps traced back
to the integrability of the recursion relations we have obtained. A first indication of this
integrability is the existence of ``integrals of motion" for these equations.
For instance, the recursion relation of the tetravalent case, eq.\quaca, has the following
integral of motion:
\eqn\intmoqua{\eqalign{ f(R_n,R_{n+1})&={\rm const.} \cr
f(x,y)&=x y (1-g x-g y)-x-y \cr}} 
as is immediately checked by forming 
\eqn\difvar{f(R_n,R_{n+1})-f(R_n,R_{n-1})=(R_{n+1}-R_{n-1})(R_n-1-g R_n(R_{n+1}+R_n+R_{n-1}))}
More generally, one may construct such integrals of motion for all the models studied in this note.
Again, it may be that this integrability property relates to the existence, for the
same planar graphs but with arbitrary topology rather than planar with fixed geodesic distances, of matrix
model formulations that are solvable via orthogonal polynomial techniques.

\noindent{\bf Acknowledgments}
This note summarizes work essentially done in collaboration with J. Bouttier and E. Guitter.
We also thank M. Bousquet-M\'elou and G. Schaeffer for fruitful discussions on constellations
and bipartite graphs.

\listrefs
\end